\documentclass[12pt,leqno]{article}

\usepackage{amsmath}
\usepackage{amssymb}
\usepackage{amsfonts}
\usepackage{amscd}
\usepackage{a4}

\allowdisplaybreaks
\numberwithin{equation}{section}
\newtheorem{Theorem}{Theorem}[section]
\newtheorem{Corollary}[Theorem]{Corollary}
\newtheorem{Lemma}[Theorem]{Lemma}
\newtheorem{Proposition}[Theorem]{Proposition}
\newtheorem{Tableau}[Theorem]{Table}

\newenvironment{Remark}{\refstepcounter{Theorem}
\noindent {\it Remark \arabic{section}.\arabic{Theorem}.}}\noindent

\newenvironment{Examples}{\refstepcounter{Theorem}
\noindent {\it Examples \arabic{section}.\arabic{Theorem}.}}\noindent

\newcommand*{\NN}{{\mathbb N}}
\newcommand*{\RR}{{\mathbb R}}
\newcommand*{\ZZ}{{\mathbb Z}}

\newcommand*{\A}{{\mathcal A}}
\newcommand*{\B}{{\mathcal B}}
\newcommand*{\F}{{\mathcal F}}
\newcommand*{\G}{{\mathcal G}}
\newcommand*{\HH}{{\mathcal H}}
\newcommand*{\I}{{\mathcal I}}
\newcommand*{\LL}{{\mathcal L}}
\newcommand*{\X}{{\mathcal X}}
\newcommand*{\Y}{{\mathcal Y}}

\newcommand*{\U}{\overline{U}}

\newcommand*{\Var}{\text{\rm Var}}

\newcommand*{\II}{I\hspace{-2pt}I}
\newcommand*{\III}{I\hspace{-2pt}I\hspace{-1.9pt}I}
\newcommand*{\IV}{I\hspace{-1.5pt}V}

\begin{document}

\renewcommand{\thefootnote}{}

\begin{center}\Large\bf 

On stable central limit theorems\\

for multivariate discrete-time martingales
\end{center}

\begin{center}
\it

Erich H\"ausler\hspace{1.5pt}$^1$\footnote{\hspace{-17.5pt}$^1$\hspace{0.4pt}Mathematical Institute, 
University of Giessen, Giessen, Germany, erich.haeusler@t-online.de} and 
Harald Luschgy\hspace{1.5pt}$^2$\footnote{\hspace{-17.5pt}$^2$\hspace{0.5pt}FB IV, Mathematics, University of Trier, Trier, Germany, luschgy@uni-trier.de}

\end{center}\bigskip\bigskip

\centerline{\it Abstract}\bigskip

We provide a systematic approach to stable central limit theorems for $d$-dimen\-sion\-al martingale difference arrays and
martingale difference sequences. The conditions imposed are straightforward extensions of the
univariate case.\bigskip

\section{\large Introduction}
\label{Section:Introduction}

It is well-known that the convergence in distribution in many classical limit theorems in probability
theory and mathematical statistics is in fact stable in the sense of Renyi, which broadens the
possible range of applications of these limit theorems. A comprehensive account of the theory of
stable convergence together with illustrative applications has been given in \cite{HaeuslerLuschgy}. 
Chapter 6 of that book is devoted to stable convergence in central limit theorems for martingale 
difference arrays and martingale difference sequences
of real-valued random variables. It is the aim of the present paper to show that the basic results presented there
also hold true for arrays and sequences of $d$-dimensional random vectors for all $d\in\NN$ if the
sufficent conditions are generalized to the $d$-dimensional case in a straightforward way. This provides
a systematic approach to stable central limit theorems for $d$-dimensional martingale difference arrays
and $d$-dimensional martingale difference sequences. For a specific example of such a result, see
\cite{CrimaldiPratelli}, Proposition~3.1.
\smallskip

The basic notations and definitions are as follows.
Let $d\in\NN$, and for every $n\in\NN$ let $k_n\in\NN$ and let $\left(X_{nk}\right)_{1\leq k\leq k_n}$
be a sequence of $\RR^d$-valued random vectors defined on
a probability space $\left(\Omega,\F,P\right)$, where $X_{nk}=\left(X_{nk,j}\right)^T_{1\leq j\leq d}$
is always read as a column vector. Moreover, let
$\left(\F_{nk}\right)_{0\leq k\leq k_n}$ be a filtration in $\F$, i.e. the $\F_{nk}$ are sub-$\sigma$-fields
of $\F$ with
$\F_{n0}\subset\F_{n1}\subset\cdots\subset\F_{nk_n}\subset\F$. The
sequence $\left(X_{nk}\right)_{1\leq k\leq k_n}$ is called
\textit{adapted to the filtration} $\left(\F_{nk}\right)_{0\leq k\leq k_n}$ if
$X_{nk}$ is measurable w.r.t. $\F_{nk}$ for all $1\leq k\leq k_n$. The triangular 
array\index{triangular array}
$\left(X_{nk}\right)_{1\leq k\leq k_n,n\in\NN}$ of $d$-dimensional random vectors is called \textit{adapted to
the triangular array}
$\left(\F_{nk}\right)_{0\leq k\leq k_n,n\in\NN}$ of $\sigma$-fields if the
row $\left(X_{nk}\right)_{1\leq k\leq k_n}$ is adapted to the filtration
$\left(\F_{nk}\right)_{0\leq k\leq k_n}$
for every $n\in\NN$. In this paper, all arrays $\left(X_{nk}\right)_{1\leq k\leq k_n,n\in\NN}$ are assumed to be
adapted to the given array $\left(\F_{nk}\right)_{0\leq k\leq k_n,n\in\NN}$ of $\sigma$-fields.
\smallskip

An array $\left(X_{nk}\right)_{1\leq k\leq k_n,n\in\NN}$ adapted to
$\left(\F_{nk}\right)_{0\leq k\leq k_n,n\in\NN}$ is called a 
\textit{martingale difference array} if $X_{nk,j}\in\LL^1\left(P\right)$  with 
$E\left(X_{nk,j}|\F_{n,k-1}\right)=0$ for all
$1\leq k\leq k_n$, $n\in\NN$ and $1\leq j\leq d$. This means that for every $n\in\NN$ 
and $1\leq j\leq d$ the sequence
$\left(X_{nk,j}\right)_{1\leq k\leq k_n}$ of real-valued random variables is a 
martingale difference sequence w.r.t. the filtration $\left(\F_{nk}\right)_{0\leq k\leq k_n}$. 
The conditional expectation of a $d$-dimensional random vector $Y=\left(Y_j\right)_{1\leq j\leq d}^T$
with $Y_j\in\LL^1\left(P\right)$ for $1\leq j\leq d$ w.r.t a sub-$\sigma$-field $\HH$ of $\F$ is defined
by $E\left(Y|\HH\right)=\left(E\left(Y_j|\HH\right)\right)_{1\leq j\leq d}^T$ so that for a martingale
difference array $\left(X_{nk}\right)_{1\leq k\leq k_n,n\in\NN}$ w.r.t. the array 
$\left(\F_{nk}\right)_{0\leq k\leq k_n,n\in\NN}$ we have $E\left(X_{nk}|\F_{n,k-1}\right)=0$
for all $n\in\NN$ and $1\leq k\leq k_n$ where $0$ is the null vector in $\RR^d$. 
Expectations and conditional expectations for random matrices are defined analogously.
A martingale difference
array is called \textit{square integrable} if all random vectors have components in $\LL^2\left(P\right)$.
\smallskip

The paper is organized as follows. In Section~\ref{Section:Threemainresults} we will present the main
results about $d$-dimensional martingale difference arrays. Section~\ref{Section:martingales} contains important consequences
for $d$-dimensional martingale difference sequences. In Section~\ref{Section:Relationship} we will formulate several
propositions concerning the relationship between different sets of conditions before giving their proofs
in Section~\ref{Section:Proofs}. In Section~\ref{Section:Application} we will investigate an application to autoregressions
of order $d$. For the readers convenience and to make this work self-contained some basic facts on
stationary and ergodic random processes are provided in the appendix.

\section{\large Main results}
\label{Section:Threemainresults}

From now on, we assume that the
sequence $\left(k_n\right)_{n\in\NN}$ is nondecreasing with $k_n\geq n$ for all $n\in\NN$.
Let $\G_{nk}=\bigcap_{m\geq n}\F_{mk}$ for $n\in\NN$ and $0\leq k\leq k_n$, and
$\G=\sigma\left(\bigcup_{n=1}^\infty\G_{nk_n}\right)$.
We always set $\F_\infty=\sigma\left(\bigcup_{n=1}^\infty\F_{nk_n}\right)$.
Let $\left\langle\cdot,\cdot\right\rangle$ denote the standard Euclidean scalar product in $\RR^d$ and let
$\left\lVert\cdot\right\rVert$ denote the associated Euclidean $\ell^2$-norm in $\RR^d$. The multivariate 
centered normal distribution with covariance matrix $\Sigma$ is denoted by 
$N\left(0,\Sigma\right)$ and $N\left(0,I_d\right)$ denotes the $d$-dimensional 
standard normal distribution. The distribution of an arbitrary random vector $X$ is denoted by
$P^X$, and equality in distribution of random vectors $X$ and $Y$ is denoted by $X\stackrel{d}{=}Y$. 
\smallskip

In our first theorem we consider martingale difference arrays under a conditional $L_1$-version of Lindeberg's classical
condition.
\medskip

\begin{Theorem}
\label{Theorem:Raikov}
Let $\left(X_{nk}\right)_{1\leq k\leq k_n,n\in\NN}$ be a martingale difference
array of $\RR^d$-valued random vectors adapted to the array $\left(\F_{nk}\right)_{0\leq k\leq k_n,n\in\NN}$.
Assume the conditional L$_1$-version of Lindeberg's condition
\begin{align}
&\sum_{k=1}^{k_n}E\left(\left\lVert X_{nk}\right\rVert1_{\left\{\left\lVert X_{nk}\right\rVert\geq\varepsilon\right\}}|\F_{n,k-1}\right)
\rightarrow0\quad\text{in probability as }n\to\infty\tag{CLB$_1^d$}\label{cond:CLB1}\\
&\qquad\text{for every }\varepsilon>0\notag
\end{align}
and Raikov's condition
\begin{align}
&\sum_{k=1}^{k_n}X_{nk}X_{nk}^T\rightarrow
A\quad\text{in probability as }n\to\infty\text{ for some random}\tag{R$^d$}\label{cond:R}\\
&\qquad\G\text{-measurable symmetric and positive semi-definite }d\times d\text{-matrix }A.\notag
\end{align}
Then
\begin{displaymath}
\sum_{k=1}^{k_n}X_{nk}\rightarrow A^{1/2}N_d\quad\G\text{-stably as }n\to\infty\,,
\end{displaymath}
where $P^{N_d}=N\left(0,I_d\right)$ and $N_d$ is independent of $\G$.
\end{Theorem}
\medskip

\begin{Remark}
\label{Remark:squareroot}
Since the map $B\rightarrow B^{1/2}$ is continuous on the set of symmetric, positive
semi-definite $d\times d$-matrices, the (uniquely determined and positive semi-definite)
square root $A^{1/2}$ of $A$ is also $\G$-measurable in the above setting. Under the
{\it nesting condition}, namely, $\F_{nk}\subset\F_{n+1,k}$ for all $n\in\NN$ and
$0\leq k\leq k_n$, we obtain $\G_{nk}=\bigcap_{m\geq n}\F_{mk}=\F_{nk}$ and thus
$\G=\F_\infty$. In this case, $\F_\infty$-measurability,
symmetry and positive semi-definiteness of $A$ (everywhere) can be assumed w.l.o.g.
\end{Remark}
\medskip

\textit{Proof of Theorem \ref{Theorem:Raikov}.} 
For $d=1$ Theorem \ref{Theorem:Raikov} follows from Corollary~6.19 and Corollary~6.22 in \cite{HaeuslerLuschgy}. 
To obtain Theorem \ref{Theorem:Raikov} for all $d\in\NN$, we will apply the
Cram\'{e}r-Wold device for stable convergence; see \cite{HaeuslerLuschgy}, Corollary 3.19, 
(i) $\Leftrightarrow$ (iii).
According to this tool for deriving the $d$-dimensional case from the one-dimensional case, we have to show that
for any $d\in\NN\setminus\left\{1\right\}$
\begin{equation}
\label{eq:convergenceofscalarproducts}
\left\langle u,\sum_{k=1}^{k_n}X_{nk}\right\rangle\rightarrow\left\langle u,A^{1/2}N_d\right\rangle\quad\G
\text{-stably as }n\to\infty\,,
\end{equation}
where $u\in\RR^d\setminus\left\{0\right\}$ is an arbitrary non-random vector. For this, let 
$u\in\RR^d\setminus\left\{0\right\}$ be fixed. For all $n\in\NN$ we have
\begin{equation*}
\left\langle u,\sum_{k=1}^{k_n}X_{nk}\right\rangle=\sum_{k=1}^{k_n}\left\langle u,X_{nk}\right\rangle=\sum_{k=1}^{k_n}Y_{nk}
\end{equation*}
with $Y_{nk}:=\left\langle u,X_{nk}\right\rangle$. Then $\left(Y_{nk}\right)_{1\leq k\leq k_n,n\in\NN}$ is a martingale
difference array of real-valued random variables adapted to $\left(\F_{nk}\right)_{0\leq k\leq k_n,n\in\NN}$.
Because of $\left\lvert Y_{nk}\right\rvert=\left\lvert\left\langle u,X_{nk}\right\rangle\right\rvert
\leq\left\lVert u\right\rVert\left\lVert X_{nk}\right\rVert$
we have for all $n\in\NN$ und $\varepsilon>0$
\begin{equation*}
\sum_{k=1}^{k_n}E\left(\left|Y_{nk}\right|1_{\left\{\left|Y_{nk}\right|\geq\varepsilon\right\}}|
\F_{n,k-1}\right)
\leq\left\lVert u\right\rVert
\sum_{k=1}^{k_n}E\left(\left\lVert X_{nk}\right\rVert1_{\left\{\left\lVert X_{nk}\right\rVert\geq\varepsilon/\left\lVert u\right\rVert\right\}}|
\F_{n,k-1}\right)\,.
\end{equation*} 
Therefore, by (\ref{cond:CLB1}),
\begin{equation*}
\sum_{k=1}^{k_n}E\left(\left|Y_{nk}\right|1_{\left\{\left|Y_{nk}\right|\geq\varepsilon\right\}}|
\F_{n,k-1}\right)\rightarrow 0\quad\text{in probability as }n\to\infty
\end{equation*}
for every $\varepsilon>0$, showing that (\ref{cond:CLB1}) holds for the martingale difference array 
$\left(Y_{nk}\right)_{1\leq k\leq k_n,n\in\NN}$ with $d=1$. Moreover, for all $n\in\NN$
\begin{displaymath}
\sum_{k=1}^{k_n}Y_{nk}^2=\sum_{k=1}^{k_n}\left\langle u,X_{nk}\right\rangle^2=\sum_{k=1}^{k_n}u^TX_{nk}X_{nk}^Tu=
u^T\left(\sum_{k=1}^{k_n}X_{nk}X_{nk}^T\right)u\rightarrow u^TAu
\end{displaymath}
in probability as $n\to\infty$
by (\ref{cond:R}). Applying now the already proven Theorem \ref{Theorem:Raikov} for $d=1$ we obtain
\begin{equation*}
\left\langle u,\sum_{k=1}^{k_n}X_{nk}\right\rangle=\sum_{k=1}^{k_n}Y_{nk}
\rightarrow\left(u^TAu\right)^{1/2}N_1\quad\G\text{-stably as }
n\to\infty\,,
\end{equation*}
where $N_1$ is a one-dimensional random variable with $P^{N_1}=N\left(0,1\right)$ 
which is independent of $\G$. According
to the definition of stable converence, see \cite{HaeuslerLuschgy}, Definition 3.15, to complete the proof of
(\ref{eq:convergenceofscalarproducts}) it remains to show that
\begin{equation}
\label{eq:distributionalequality}
P^{\left(u^TAu\right)^{1/2}N_1|\G}=P^{\left\langle u,A^{1/2}N_d\right\rangle|\G}\,.
\end{equation}
For a deterministic matrix $A$ 
we have $P^{\left(u^TAu\right)^{1/2}N_1|\G}=N\left(0,u^TAu\right)$
where
\begin{equation*}
u^TAu=u^TA^{1/2}A^{1/2}u=\left(A^{1/2}u\right)^TA^{1/2}u=\left\langle A^{1/2}u,A^{1/2}u\right\rangle
=\left\lVert A^{1/2}u\right\rVert^2.
\end{equation*}
Hence
\begin{equation*}
\left\langle u,A^{1/2}N_d\right\rangle=\left\langle A^{1/2}u,N_d\right\rangle\stackrel{d}{=}
\left\lVert A^{1/2}u\right\rVert N_1=\left(u^TAu\right)^{1/2}N_1\,.
\end{equation*}
For a random $A$, both $A$ and $A^{1/2}$ are $\G$-measurable so that 
\eqref{eq:distributionalequality} holds as well.
Thus, Theorem \ref{Theorem:Raikov} is proven.\hfill$\Box$
\medskip

The most common and useful multivariate stable central limit theorem stated in the following
theorem is a consequence of Theorem \ref{Theorem:Raikov}.
\medskip

\begin{Theorem}
\label{Theorem:Lindeberg}
In the situation of Theorem \ref{Theorem:Raikov} let the martingale difference array
$\left(X_{nk}\right)_{1\leq k\leq k_n,n\in\NN}$ be square integrable. Assume the
conditional form of Lindeberg's condition
\begin{align}
&\sum_{k=1}^{k_n}E\left(\left\lVert X_{nk}\right\rVert^2
1_{\left\{\left\lVert X_{nk}\right\rVert\geq\varepsilon\right\}}|\F_{n,k-1}\right)
\rightarrow0\quad\text{in probability as }n\to\infty\tag{CLB$^d$}\label{cond:CLB}\\
&\qquad\text{for every }\varepsilon>0\notag
\end{align}
and
\begin{align}
&\sum_{k=1}^{k_n}E\left(X_{nk}X_{nk}^T|\F_{n,k-1}\right)\rightarrow
A\quad\text{in probability as }n\to\infty\text{ for some random}\tag{N$^d$}\label{cond:N}\\
&\qquad\G\text{-measurable symmetric and positive semi-definite }d\times d\text{-matrix }A.\notag
\end{align}
Then
\begin{displaymath}
\sum_{k=1}^{k_n}X_{nk}\rightarrow A^{1/2}N_d\quad\G\text{-stably as }n\to\infty\,,
\end{displaymath}
where $N_d$ is as in Theorem~\ref{Theorem:Raikov}.
\end{Theorem}
\medskip

\textit{Proof.} By the subsequent Proposition \ref{Proposition:CLBandNimplyR}, 
the conditions \eqref{cond:CLB} and
\eqref{cond:N} imply the conditions \eqref{cond:CLB1} and \eqref{cond:R}. Consequently,
Theorem~\ref{Theorem:Lindeberg} follows from Theorem~\ref{Theorem:Raikov}.
\smallskip

Alternatively, one may argue along the lines of the proof of Theorem~\ref{Theorem:Raikov}
based on the corresponding univariate stable central limit theorem (see
\cite{HaeuslerLuschgy}, Theorem~6.1).

\hfill$\Box$
\medskip

Our final result for arrays $\left(X_{nk}\right)_{1\leq k\leq k_n,n\in\NN}$ 
holds for arbitrary adapted arrays without any moment assumptions and provides a 
generalization of Theorem~\ref{Theorem:Raikov}.
For $0<a<\infty$ and $1\leq k\leq k_n,n\in\NN$ we set
\begin{equation}
\label{eq:Xnktruncateda}
X_{nk}\left(a\right)=X_{nk}1_{\left\{\left\lVert X_{nk}\right\rVert\leq a\right\}}-
E\left(X_{nk}1_{\left\{\left\lVert X_{nk}\right\rVert\leq a\right\}}
|\F_{n,k-1}\right)\,.
\end{equation}
\medskip

\begin{Theorem}
\label{Theorem:General}
Let the array $\left(X_{nk}\right)_{1\leq k\leq k_n,n\in\NN}$ 
be adapted to 
$\left(\F_{nk}\right)_{0\leq k\leq k_n,n\in\NN}$. Assume that for some $a>0$
\begin{align}
&\sum_{k=1}^{k_n}X_{nk}1_{\left\{\left\lVert X_{nk}\right\rVert>a\right\}}+
E\left(X_{nk}1_{\left\{\left\lVert X_{nk}\right\rVert\leq a\right\}}|\F_{n,k-1}\right)\rightarrow0
\tag{T$_a^d$}\label{cond:Ta}\\
&\qquad\text{in probability as }n\to\infty\,,\notag
\end{align}
\begin{align}
&\max_{1\leq k\leq k_n}\left\lVert X_{nk}\left(a\right)\right\rVert\rightarrow0\quad\text{in probability as }n\to\infty
\tag{TM$_a^d$}\label{cond:TMa}
\end{align}
and
\begin{align}
&\sum_{k=1}^{k_n}X_{nk}\left(a\right)X_{nk}\left(a\right)^T\rightarrow A\quad\text{in probability as }n\to\infty
\text{ for some }\tag{TR$_a^d$}\label{cond:TRa}\\
&\qquad\text{random }\G\text{-measurable symmetric and positive semi-definite}\notag\\
&\qquad d\times d\text{-matrix }A\,.\notag
\end{align}
Then
\begin{displaymath}
\sum_{k=1}^{k_n}X_{nk}\rightarrow A^{1/2}N_d\quad\G\text{-stably as }n\to\infty\,,
\end{displaymath}
where $N_d$ is as in Theorem \ref{Theorem:Raikov}.
\end{Theorem}
\medskip

\begin{Remark}
\label{Remark:general}
The subsequent Proposition \ref{Proposition:CLB1andRimplytheT} shows that 
for martingale difference arrays the conditions
\eqref{cond:CLB1} and \eqref{cond:R} imply the conditions \eqref{cond:Ta}, \eqref{cond:TMa}
and \eqref{cond:TRa} for every $a>0$. Hence, Theorem~\ref{Theorem:Raikov} follows in fact from 
Theorem~\ref{Theorem:General}.
\end{Remark}
\medskip

\textit{Proof of Theorem \ref{Theorem:General}}. Let $\left(X_{nk}\right)_{1\leq k\leq k_n,n\in\NN}$ be an
array of $d$-dimensional random vectors satisfying \eqref{cond:Ta}, \eqref{cond:TMa} and \eqref{cond:TRa} for
some $0<a<\infty$, and let $\left(X_{nk}\left(a\right)\right)_{1\leq k\leq k_n,n\in\NN}$ 
be the bounded martingale difference array
defined by \eqref{eq:Xnktruncateda}. By Proposition \ref{Proposition:Approx} {\rm(a)} 
the latter array satisfies conditions \eqref{cond:CLB1}
and \eqref{cond:R}. Consequently, Theorem \ref{Theorem:Raikov} gives
\begin{displaymath}
\sum_{k=1}^{k_n}X_{nk}\left(a\right)\rightarrow A^{1/2}N_d\quad\G\text{-stably as }n\to\infty\,.
\end{displaymath}
Because \eqref{eq:approxXnka} is true by 
Proposition \ref{Proposition:Approx}~(a), an application of Theorem~3.18~(a) in \cite{HaeuslerLuschgy} yields
\begin{displaymath}
\sum_{k=1}^{k_n}X_{nk}\rightarrow A^{1/2}N_d\quad\G\text{-stably as }n\to\infty\,.
\end{displaymath}
This completes the proof of Theorem~\ref{Theorem:General}.\hfill$\Box$
\medskip

\section{Martingales}
\label{Section:martingales}

Now we consider applications to martingales. For this, let
$\left(X_k\right)_{k\geq1}$ be an infinite sequence of integrable $d$-dimensional random vectors
$X_k=\left(X_{k,j}\right)^T_{1\leq j\leq d}$ for $d\in\NN$ which is a martingale difference sequence w.r.t.
the infinite filtration $\left(\F_k\right)_{k\geq0}$. Setting $\F_{nk}=\F_k$ for all $n\in\NN$ and
$0\leq k\leq n$ we obtain an array $\left(\F_{nk}\right)_{0\leq k\leq n,n\in\NN}$ of $\sigma$-fields
with $k_n=n$. Let
$\left(K_n\right)_{n\in\NN}$ be a sequence of non-random $d\times d$-matrices
$K_n=\left(K_{n,j\ell}\right)_{1\leq j,\ell\leq d}$. For all $n\in\NN$ and $1\leq k\leq n$ we set
$X_{nk}=K_nX_k$. Then $\left(X_{nk}\right)_{1\leq k\leq n,n\in\NN}$ is a array of $d$-dimensional
random vectors with
\begin{equation}
\label{eq:representationXnkj}
X_{nk,j}=\sum_{\ell=1}^d K_{n,j\ell}X_{k,\ell}\quad\text{for all }n\in\NN, 1\leq k\leq n\text{ and }
1\leq j\leq d\,.
\end{equation}
From \eqref{eq:representationXnkj} we see that $\left(X_{nk}\right)_{1\leq k\leq n,n\in\NN}$ is a martingale
difference array w.r.t. the array $\left(\F_{nk}\right)_{0\leq k\leq n,n\in\NN}$ of $\sigma$-fields. 
With $\F_\infty=\sigma\left(\bigcup_{k=0}^\infty\F_k\right)$ we have the following 
corollary to Theorem~\ref{Theorem:Raikov}.
\medskip

\begin{Corollary}{\rm(}Martingales{\rm)}
\label{Corollary:Raikovforsequences}
If
\begin{align*}
&\sum_{k=1}^nE\left(\left\lVert K_nX_k\right\rVert1_{\{\left\lVert K_nX_k\right\rVert\geq\varepsilon\}}|\F_{k-1}\right)\rightarrow0
\quad\text{in probability as }n\to\infty\tag{CLB$_{1,K_n}^d$}\label{eq:CLB1forsequences}\\
&\qquad \text{for every }\varepsilon>0\notag
\end{align*}
and
\begin{align*}
K_n\sum_{k=1}^nX_kX_k^TK_n^T\rightarrow A\quad\text{in probability as }n\to\infty
\text{ for some random}\quad\,\,\tag{R$_{K_n}^d$}\label{eq:Raikovforsequences}\\
\qquad\F_\infty\text{-measurable symmetric and positive semi-definite }d\times d\text{-matrix }A\,,
\end{align*}
then
\begin{displaymath}
K_n\sum_{k=1}^nX_k\rightarrow A^{1/2}N_d\quad\F_\infty\text{-stably as }n\to\infty\,,
\end{displaymath}
where $P^{N_d}=N\left(0,I_d\right)$ and $N_d$ is independent of $\F_\infty$.
\end{Corollary}
\medskip

Note that in the martingale setting the nesting condition holds so that we have 
$\G_{nk}=\bigcap_{m\geq n}\F_{mk}=\bigcap_{m\geq n}\F_k=\F_k$ and
$\G=\sigma\left(\bigcup_{n=1}^\infty\G_{nn}\right)=\sigma\left(\bigcup_{n=1}^\infty\F_n\right)=\F_\infty$.
Condition~\eqref{eq:CLB1forsequences} clearly is condition~\eqref{cond:CLB1} 
for the array $\left(X_{nk}\right)_{1\leq k\leq n,n\in\NN}$,
and~\eqref{eq:Raikovforsequences} is condition~\eqref{cond:R} with $\G=\F_\infty$. Therefore, the corollary is 
a special case of Theorem~\ref{Theorem:Raikov}. Since the condition~\eqref{eq:CLB1forsequences}
is weaker than the condition
\begin{equation*}
E\left(\max_{1\leq k\leq n}\left\lVert K_nX_k\right\rVert\right)\rightarrow0\quad\text{as }n\to\infty
\tag{M$_{1,K_n}^d$}
\end{equation*}
(see Proposition~\ref{Proposition:MimpliesCLB}), Corollary~\ref{Corollary:Raikovforsequences}
provides a discrete-time version of Theorem~2.2 in \cite{CrimaldiPratelli}.
\smallskip

Theorem~\ref{Theorem:Lindeberg} gives the following result.
\medskip

\begin{Corollary}{\rm(}Martingales{\rm)}
\label{Corollary:Lindebergforsequences}
Let $\left(X_k\right)_{k\geq1}$ be a square integrable martingale difference 
sequence w.r.t. $\left(\F_k\right)_{k\geq0}$.
If
\begin{align*}
&\sum_{k=1}^nE\left(\left\lVert K_nX_k\right\rVert^21_{\{\left\lVert K_nX_k\right\rVert\geq\varepsilon\}}|\F_{k-1}\right)\rightarrow0
\quad\text{in probability as }n\to\infty\tag{CLB$_{K_n}^d$}\label{eq:CLBforsequences}\\
&\qquad \text{for every }\varepsilon>0
\end{align*}
and
\begin{align*}
&K_n\sum_{k=1}^nE\left(X_kX_k^T|\F_{k-1}\right)K_n^T\rightarrow A\quad\text{in probability as }
n\to\infty\text{ for some}\tag{N$_{K_n}^d$}\label{eq:normingforsequences}\\
&\text{random }\F_\infty\text{-measurable symmetric and positive semi-definite}\\
&d\times d\text{-matrix }A\,,
\end{align*}
then
\begin{displaymath}
K_n\sum_{k=1}^nX_k\rightarrow A^{1/2}N_d\quad\F_\infty\text{-stably as }n\to\infty\,,
\end{displaymath}
where $N_d$ is as in Corollary~\ref{Corollary:Raikovforsequences}.
\end{Corollary}
\medskip

Our next result is for $d$-dimensional stationary martingale difference sequences. 
The shift on $\left(\RR^d\right)^\NN$ is defined by 
$S\left(x\right)=S\left(\left(x_n\right)_{n\geq1}\right)=\left(x_{n+1}\right)_{n\geq1}$. Note that it is
$\left(\B\left(\RR^d\right)^\NN,\B\left(\RR^d\right)^\NN\right)$-measurable, where $\B\left(\RR^d\right)$ denotes
the Borel $\sigma$-field on $\RR^d$. Let
\begin{displaymath}
\I=\B\left(\RR^d\right)^\NN\left(S\right)=\left\{A\in\B\left(\RR^d\right)^\NN:S^{-1}\left(A\right)=A\right\}
\end{displaymath}
denote the $\sigma$-field of invariant measurable subsets of $\left(\RR^d\right)^\NN$. The Frobenius norm
(or Euclidean $\ell^2$-norm) on the space $\RR^{m\times n}$ of all real $m\times n$-matrices $D$ is given by
$\left\lVert D\right\rVert_F=\left(\sum_{i=1}^m\sum_{j=1}^nD_{ij}^2\right)^{1/2}$.
\medskip

\begin{Corollary}{\rm(}Stationary martingale differences{\rm)}
\label{Corollary:Stationarymartingaledifferences}
Let $X=\left(X_n\right)_{n\geq1}$ be an $\RR^d$-valued stationary martingale difference sequence w.r.t. 
$\left(\F_n\right)_{n\geq0}$ with $\lVert X_1\rVert\in\LL^2\left(P\right)$. Then
\begin{displaymath}
\frac{1}{\sqrt{n}}\sum_{k=1}^nX_k\rightarrow E\left(X_1X_1^T|\I_X\right)^{1/2}N_d\quad\F_\infty\text{-stably as }n\to\infty\,,
\end{displaymath}
where $\I_X=X^{-1}\left(\I\right)$, $N_d$ is independent of $\F_\infty$ and $P^{N_d}=N\left(0,I_d\right)$. 
If $X$ is additionally ergodic, then
\begin{displaymath}
\frac{1}{\sqrt{n}}\sum_{k=1}^nX_k\rightarrow E\left(X_1X_1^T\right)^{1/2}N_d\quad\F_\infty\text{-mixing as }n\to\infty\,.
\end{displaymath}
\end{Corollary}
\medskip

\textit{Proof.} We apply Corollary~\ref{Corollary:Raikovforsequences} with $K_n=n^{-1/2}I_d$. By stationarity
of $X$
\begin{align*}
&E\left(\frac{1}{n}\sum_{k=1}^nE\left(\lVert X_k\rVert^21_{\left\{\lVert X_k\rVert\geq\varepsilon n^{1/2}\right\}}|\F_{k-1}\right)\right)
=\frac{1}{n}\sum_{k=1}^nE\left(\left\lVert X_k\right\rVert^2
1_{\left\{\left\lVert X_k\right\rVert\geq\varepsilon n^{1/2}\right\}}\right)\\
&\qquad=E\left(\lVert X_1\rVert^21_{\left\{\lVert X_1\rVert\geq\varepsilon n^{1/2}\right\}}\right)\rightarrow0
\end{align*}
as $n\to\infty$
because of $E\left(\lVert X_1\rVert^2\right)<\infty$ so that condition \eqref{eq:CLBforsequences} is satisfied.
By Proposition~\ref{Proposition:CLBandNimplyR}, this condition implies \eqref{eq:CLB1forsequences}.
As for condition \eqref{eq:Raikovforsequences}, we apply the ergodic theorem (see Theorem~\ref{App:Birkhoff} with
$\left(\X,\A\right)=\left(\RR^d,\B\left(\RR^d\right)\right)$, $T=\NN$ and
$f:\left(\RR^d\right)^\NN\rightarrow\RR^{d\times d},\ f\left(x\right)=\pi_1\left(x\right)\pi_1\left(x\right)^T$,
where $\pi_1:\left(\RR^d\right)^\NN\rightarrow\RR^d,\ \pi_1\left(x\right)=x_1$, so that
$\lVert f\left(X\right)\rVert_F\leq\lVert X_1\rVert^2$). This yields
\begin{displaymath}
\frac{1}{n}\sum_{k=1}^nX_kX_k^T\rightarrow E\left(X_1X_1^T|\I_X\right)\quad\text{almost surely as }n\to\infty\,.
\end{displaymath}
Now the assertion follows from Corollary~\ref{Corollary:Raikovforsequences} and $E\left(X_1X_1^T|\I_X\right)=
E\left(X_1X_1^T\right)$ almost surely in the ergodic case.\hfill$\Box$
\medskip

In our final result in this section we consider a sequence $\left(\HH_k\right)_{k\geq1}$ of sub-$\sigma$-fields
of $\F$ with $\HH_k\supset\HH_{k+1}$ for every $k\geq1$, which we call a \textit{reversed filtration}.
Set $\HH_\infty:=\bigcap_{k\geq1}\HH_k$. We call a sequence $\left(X_k\right)_{k\geq1}$ of integrable 
$d$-dimensional random vectors which is adapted to $\left(\HH_k\right)_{k\geq1}$ (i.e. $X_k$ is $\HH_k$-measurable
for all $k\geq1$) a \textit{reversed martingale difference sequence} if $E\left(X_k|\HH_{k+1}\right)=0$ for all
$k\geq1$. Note that in general the partial sum $\sum_{k=1}^nX_k$ will not be $\HH_n$-measurable which means that
the sequence $\left(\sum_{k=1}^nX_k\right)_{n\geq1}$ will not be a reversed martingale w.r.t. $\left(\HH_k\right)_{k\geq1}$.
Therefore, the following limit theorem should not be confused with central limit theorems for reversed martingales
as obtained e.g. in \cite{Loynes}.
\medskip

\begin{Corollary}{\rm(}reversed martingale differences{\rm)}
\label{Corollary:Reversedmartingaledifferences}
Let $\left(X_k\right)_{k\geq1}$ be a square-inte\-grable $\RR^d$-valued reversed martingale difference sequence w.r.t. 
the reversed filtration
$\left(\HH_k\right)_{k\geq1}$. If
\begin{align*}
&\sum_{k=1}^nE\left(\left\lVert K_nX_k\right\rVert^21_{\{\left\lVert K_nX_k\right\rVert\geq\varepsilon\}}|\HH_{k+1}\right)\rightarrow0
\quad\text{in probability as }n\to\infty\tag{CLB$_{K_n}^{d,r}$}\label{cond:CLBforinversesequences}\\
&\qquad \text{for every }\varepsilon>0
\end{align*}
and
\begin{align*}
&K_n\sum_{k=1}^nE\left(X_kX_k^T|\HH_{k+1}\right)K_n^T\rightarrow A\quad\text{in probability as }
n\to\infty\text{ for some}\tag{N$_{K_n}^{d,r}$}\label{cond:normingforinversesequences}\\
&\text{random }\HH_\infty\text{-measurable symmetric and positive semi-definite}\\
&d\times d\text{-matrix }A\,,
\end{align*}
then
\begin{displaymath}
K_n\sum_{k=1}^nX_k\rightarrow A^{1/2}N_d\quad\HH_\infty\text{-stably as }n\to\infty\,,
\end{displaymath}
where $P^{N_d}=N\left(0,I_d\right)$ and $N_d$ is independent of $\HH_\infty$.
\end{Corollary}
\medskip

\textit{Proof.} For all $n\in\NN$ and $0\leq k\leq n$ set $\F_{nk}:=\HH_{n+1-k}$. Then 
$\left(\F_{nk}\right)_{0\leq k\leq n,n\in\NN}$ is a triangular array of $\sigma$-fields as defined in the
introduction. Let $X_{nk}=K_nX_{n+1-k}$ for all $n\in\NN$ and $1\leq k\leq n$. Then 
$\left(X_{nk}\right)_{1\leq k\leq n,n\in\NN}$ is a square integrable martingale difference array
w.r.t. $\left(\F_{nk}\right)_{0\leq k\leq n,n\in\NN}$. For all $n\in\NN$ and $\varepsilon>0$ we have
\begin{align*}
&\sum_{k=1}^nE\left(\left\lVert X_{nk}\right\rVert^21_{\{\left\lVert X_{nk}\right\rVert\geq\varepsilon\}}|\F_{k-1}\right)=
\sum_{k=1}^nE\left(\left\lVert K_nX_{n+1-k}\right\rVert^21_{\{\left\lVert K_nX_{n+1-k}\right\rVert\geq\varepsilon\}}|\HH_{n+2-k}\right)\\
&\quad=\sum_{k=1}^nE\left(\left\lVert K_nX_k\right\rVert^21_{\{\left\lVert K_nX_k\right\rVert\geq\varepsilon\}}|\HH_{k+1}\right)
\end{align*}
which converges to $0$ in probability as $n\to\infty$ by \eqref{cond:CLBforinversesequences} so that
\eqref{cond:CLB} is satisfied for the array $\left(X_{nk},\F_{nk}\right)$. Moreover, for all $n\in\NN$
\begin{align*}
&\sum_{k=1}^nE\left(X_{nk}X_{nk}^T|\F_{n,k-1}\right)=K_n\sum_{k=1}^nE\left(X_{n+1-k}X_{n+1-k}^T|\HH_{n+2-k}\right)K_n^T\\
&\quad=K_n\sum_{k=1}^nE\left(X_kX_k^T|\HH_{k+1}\right)K_n^T
\end{align*}
which converges to $A$ in probability as $n\to\infty$ by \eqref{cond:normingforinversesequences} so that \eqref{cond:N}
is satisfied. In the present setting we have $\G_{nk}=\bigcap_{m\geq n}\F_{mk}=\bigcap_{m\geq n}\HH_{m+1-k}=\HH_\infty$
for all $n\in\NN$ and $0\leq k\leq n$ so that $\G=\sigma\left(\bigcup_{n=1}^\infty\G_{nn}\right)=\HH_\infty$.
Therefore, Theorem~\ref{Theorem:Lindeberg} implies
\begin{align*}
\sum_{k=1}^nX_{nk}\rightarrow A^{1/2}N_d\quad\HH_\infty\text{-stably as }n\to\infty
\end{align*}
which in view of
\begin{align*}
\sum_{k=1}^nX_{nk}=K_n\sum_{k=1}^nX_{n+1-k}=K_n\sum_{k=1}^nX_k
\end{align*}
is the assertion of the corollary.\hfill$\Box$
\medskip

Clearly, Corollary~\ref{Corollary:Reversedmartingaledifferences} is the analogue of Corollary~\ref{Corollary:Lindebergforsequences}
for reversed martingale difference sequences. An analogue of Corollary~\ref{Corollary:Raikovforsequences} for reversed
martingale difference sequences can be obtained in the same way, which in turn implies an analogue of 
Corollary~\ref{Corollary:Stationarymartingaledifferences} for stationary reversed martingale difference sequences.
\medskip

A special case of Corollary~\ref{Corollary:Reversedmartingaledifferences} for $d=1$ with convergence in distribution
instead of stable convergence has been applied in \cite{ConzeRaugi}; see Theorem~5.8 of that paper. The reference to
``reversed martingales" in the formulation of that theorem is somewhat irritating because reversed martingales play no
role there, in the same way as they play no role in our Corollary~\ref{Corollary:Reversedmartingaledifferences}.
\medskip

\section{\large Sufficient conditions}
\label{Section:Relationship}

In this section we will disclose the relationship between the sufficient conditions in 
Theorems~\ref{Theorem:Raikov}, \ref{Theorem:Lindeberg} and
\ref{Theorem:General} as well as the additional condition
\begin{align}
E\left(\max_{1\leq k\leq k_n}\left\lVert X_{nk}\right\rVert^2\right)\rightarrow0\quad
\text{as }n\to\infty\tag{M$_2^d$}\label{cond:M2}
\end{align}
for square integrable arrays $\left(X_{nk}\right)_{1\leq k\leq k_n,n\in\NN}$ and the condition
\begin{align}
E\left(\max_{1\leq k\leq k_n}\left\lVert X_{nk}\right\rVert\right)\rightarrow0\quad
\text{as }n\to\infty\tag{M$_1^d$}\label{cond:M1}
\end{align}
for integrable arrays. It will turn out that the implications contained in the following table are true.
\medskip

\begin{Tableau} 
\label{Table}
\textit{ Relationship between conditions in Theorems~\ref{Theorem:Raikov},
\ref{Theorem:Lindeberg} and \ref{Theorem:General}}
\begin{center}
\begin{tabular}{ccccc}
\eqref{cond:M2} and \eqref{cond:N} & $\stackrel{(a)}{\Rightarrow}$ & \eqref{cond:M1} and \eqref{cond:R} & & \\[8pt]
$\Downarrow{\scriptstyle(b)}$ & & $\Downarrow{\scriptstyle(c)}$ & & \\[5pt]
\eqref{cond:CLB} and \eqref{cond:N} & $\stackrel{(d)}{\Rightarrow}$ & \eqref{cond:CLB1} and 
\eqref{cond:R} & $\stackrel{(e)}{\Rightarrow}$ & 
\eqref{cond:Ta}, \eqref{cond:TMa} and \eqref{cond:TRa}\\[8pt]
& & & & for all $0<a<\infty$
\end{tabular}
\end{center}
\end{Tableau}
\medskip

The conditions in the left column require square integrable random vectors, in the middle 
integrability is sufficient, and on
the right-hand side no moment conditions are needed at all. The array 
$\left(X_{nk}\right)_{1\leq k\leq k_n,n\in\NN}$ is assumed
to be adapted to the array $\left(\F_{nk}\right)_{0\leq k\leq k_n,n\in\NN}$, and for implication 
$\left(e\right)$ it has to 
be a martingale difference array. The implications in Table~\ref{Table} are proven in a sequence of propositions.
The proofs are given in Section \ref{Section:Proofs}.
\medskip

\begin{Proposition}
\label{Proposition:MimpliesCLB}
For square integrable adapted arrays we have {\rm(\ref{cond:M2})}\,$\Rightarrow$\,{\rm(\ref{cond:CLB})}, 
and for integrable arrays we have
{\rm(\ref{cond:M1})}$\,\Rightarrow\,${\rm(\ref{cond:CLB1})}.
\end{Proposition}
\medskip

Proposition \ref{Proposition:MimpliesCLB} shows that the implications $(b)$ and $(c)$ in Table~\ref{Table} are true.
The next result is a crucial tool in the proofs of Propositions~\ref{Proposition:CLBandNimplyR} 
and \ref{Proposition:CLBandRimplyN}.
\medskip

\begin{Lemma}
\label{Lemma}
Let $d=2$. Let $\left(X_{nk}\right)_{1\leq k\leq k_n,n\in\NN}$ be a $2$-dimensional square
integrable adapted array. 
The conditions
\begin{equation}\label{eq:CLBL2Znk}
\sum_{k=1}^{k_n}E\left(\left\lVert X_{nk}\right\rVert^2
1_{\left\{\left\lVert X_{nk}\right\rVert\geq\varepsilon\right\}}\right)\rightarrow0
\quad \text{for all }\varepsilon>0
\end{equation}
and
\begin{equation}\label{eq:L2boundZnk}
\sum_{k=1}^{k_n}E\left(\left\lVert X_{nk}\right\rVert^2\right)\leq C<\infty\quad
\text{for some constant }C\text{ and all }n\in\NN
\end{equation}
imply
\begin{displaymath}
\sum_{k=1}^{k_n}E\left(X_{nk,1}X_{nk,2}|\F_{n,k-1}\right)-\sum_{k=1}^{k_n}X_{nk,1}X_{nk,2}
\rightarrow0\quad\text{in }L_1\text{ as }n\to\infty.
\end{displaymath}
\end{Lemma}
\medskip

\begin{Proposition}
\label{Proposition:CLBandNimplyR}
For square integrable adapted arrays condition \eqref{cond:CLB} implies \eqref{cond:CLB1} and conditions {\rm(\ref{cond:CLB})} and 
{\rm(\ref{cond:N})} imply {\rm(\ref{cond:R})}.
\end{Proposition}
\medskip

Clearly, implication $(d)$ in Table~\ref{Table} follows from Proposition~\ref{Proposition:CLBandNimplyR}. 
Note that from the already proven
implications $(b)$ and $(d)$ in Table~\ref{Table} we see that conditions (\ref{cond:M2}) and (\ref{cond:N}) imply
(\ref{cond:R}). Because condition (\ref{cond:M1}) is trivially implied by (\ref{cond:M2}), implication $(a)$ in
Table~\ref{Table} is also true.
\smallskip

For uniformly bounded arrays satisfying~\eqref{cond:CLB} the conditions~\eqref{cond:N} 
and~\eqref{cond:R} are in fact equivalent, as shown by the following proposition.
\medskip

\begin{Proposition}
\label{Proposition:CLBandRimplyN}
For uniformly bounded adapted arrays, that is, $\left\lVert X_{nk}\right\rVert\leq a$ for all $n\in\NN$ 
and $1\leq k\leq k_n$ and some $a<\infty$,
conditions~\eqref{cond:CLB} and~\eqref{cond:R} imply~\eqref{cond:N}.
\end{Proposition}
\medskip

Of course, in the preceding setting the conditions \eqref{cond:CLB} and \eqref{cond:CLB1} are
equivalent.
\smallskip

The next result proves implication $(e)$ in Table~\ref{Table}.
\medskip

\begin{Proposition}
\label{Proposition:CLB1andRimplytheT}
Let $\left(X_{nk}\right)_{1\leq k\leq k_n,n\in\NN}$ be a martingale difference array w.r.t the
$\sigma$-fields $\left(\F_{nk}\right)_{0\leq k\leq k_n,n\in\NN}$ which satisfies 
{\rm(\ref{cond:CLB1})} and {\rm(\ref{cond:R})}. Then
{\rm(\ref{cond:Ta})}, {\rm(\ref{cond:TMa})} and {\rm(\ref{cond:TRa})} are true for all $a>0$ 
and the limiting matrix $A$ in {\rm(\ref{cond:R})}
and {\rm(\ref{cond:TRa})} is the same.
\end{Proposition} 
\medskip

Clearly, Proposition~\ref{Proposition:CLB1andRimplytheT} shows that implication $\left(e\right)$ 
in Table~\ref{Table} is true, and all
implications in Table~\ref{Table} are proven. The first part of the next proposition will be crucial for 
deriving Theorem~\ref{Theorem:General} from Theorem~\ref{Theorem:Raikov} as well as from
Theorem~\ref{Theorem:Lindeberg}.
\medskip

\begin{Proposition}
\label{Proposition:Approx}
Let $\left(X_{nk}\right)_{1\leq k\leq k_n,n\in\NN}$ be adapted to the $\sigma$-fields 
$\left(\F_{nk}\right)_{0\leq k\leq k_n,n\in\NN}$.
Assume that for some $a>0$ conditions {\rm(\ref{cond:Ta})}, {\rm(\ref{cond:TMa})} and 
{\rm(\ref{cond:TRa})} are satisfied. 

{\rm(a)} The array $\left(X_{nk}\left(a\right)\right)_{1\leq k\leq k_n,n\in\NN}$ is a uniformly 
bounded martingale difference array w.r.t.
$\left(\F_{nk}\right)_{0\leq k\leq k_n,n\in\NN}$ which satisfies the conditions 
{\rm(\ref{cond:CLB1})} and {\rm(\ref{cond:R})} as well as
\begin{equation}
\label{eq:approxXnka}
\sum_{k=1}^{k_n}X_{nk}-\sum_{k=1}^{k_n}X_{nk}\left(a\right)\rightarrow0\quad\text{in probability as }n\to\infty\,.
\end{equation}

{\rm(b)} The array $\left(X_{nk}\left(a\right)\right)_{1\leq k\leq k_n,n\in\NN}$ also
satisfies conditions \eqref{cond:M2} and \eqref{cond:N}.
\end{Proposition}
\medskip

The second part of Proposition~\ref{Proposition:Approx} says that an array 
$\left(X_{nk}\right)_{1\leq k\leq k_n}$ of
$d$-dimen\-sion\-al random vectors without any integrability assumptions
which satisfies the three conditions~\eqref{cond:Ta}, \eqref{cond:TMa} and \eqref{cond:TRa} for some $0<a<\infty$
appearing in the lower right corner of Table~\ref{Table} is, as far as stable convergence of the row-sums is concerned, 
asymptotically equivalent to a uniformly bounded martingale difference array
satisfying the conditions in the upper left corner of Table~\ref{Table}, which formally is the strongest set of conditions. 
Consequently, Table~\ref{Table} shows that as soon as any one of the three Theorems~\ref{Theorem:Raikov},
\ref{Theorem:Lindeberg} or \ref{Theorem:General} is proven, the other two follow by appropriate implications 
in Table~\ref{Table}. We have chosen here to derive Theorem~\ref{Theorem:Raikov} directly from the one-dimensional 
case by an application of the Cram\'er-Wold device
and to obtain Theorems~\ref{Theorem:Lindeberg} and \ref{Theorem:General} indirectly via Table~\ref{Table}, but
other approaches are possible as well. Therefore,
in an informal sense, all conditions appearing in Table~\ref{Table} may be considered as being tantamount to each other. 
\smallskip

In our final result in this section we will show that the conditional Lindeberg condition 
\eqref{cond:CLB} holds if and only if it
holds for every component of the random vectors $X_{nk}$. This result is not needed in this 
paper but it may turn out to be useful in applications.
\medskip

\begin{Proposition}
\label{Proposition:EquivalenceforCLB}
For any array $\left(X_{nk}\right)_{1\leq k\leq k_n,n\in\NN}$ of $d$-dimensional random 
vectors the following statements are equivalent:
\begin{equation*}
\sum_{k=1}^{k_n}E\left(\left\lVert X_{nk}\right\rVert^2
1_{\{\left\lVert X_{nk}\right\rVert\geq\varepsilon\}}|\F_{n,k-1}\right)\rightarrow0\quad
\text{in probability as }n\to\infty\leqno{(\rm{i})}
\end{equation*}
\begin{displaymath}
\text{for every }\varepsilon>0\,.\hspace{200pt}
\end{displaymath}
{\rm(ii)} For every $1\leq j\leq d$
\begin{equation*}
\sum_{k=1}^{k_n}E\left(\left|X_{nk,j}\right|^2
1_{\{\left|X_{nk,j}\right|\geq\varepsilon\}}|\F_{n,k-1}\right)\rightarrow0\quad
\text{in probability as }n\to\infty
\end{equation*}
\begin{displaymath}
\text{for every }\varepsilon>0\,.\hspace{200pt}
\end{displaymath}
The same holds for conditions \eqref{cond:CLB1}, \eqref{cond:M1} and \eqref{cond:M2}. 
\end{Proposition}
\medskip

\section{\large Proofs of the results of Section~\ref{Section:Relationship}}
\label{Section:Proofs}

\textit{Proof of Proposition \ref{Proposition:MimpliesCLB}.}
For the proof of (\ref{cond:M2})\,$\Rightarrow$\,(\ref{cond:CLB}) note that for all 
$\varepsilon,\delta>0$ and $n\in\NN$
\begin{displaymath}
P\left(\sum_{k=1}^{k_n}\left\lVert X_{nk}\right\rVert^2
1_{\left\{\left\lVert X_{nk}\right\rVert\geq\varepsilon\right\}}\geq\delta\right)\leq
P\left(\max_{1\leq k\leq k_n}\left\lVert X_{nk}\right\rVert\geq\varepsilon\right)\,.
\end{displaymath}
Therefore, (\ref{cond:M2}) implies
\begin{displaymath}
\sum_{k=1}^{k_n}\left\lVert X_{nk}\right\rVert^2
1_{\left\{\left\lVert X_{nk}\right\rVert\geq\varepsilon\right\}}\rightarrow0\quad\text{in probability as }
n\to\infty\,.
\end{displaymath}
Because (\ref{cond:M2}) implies also uniform integrability of 
$\left\{\max_{1\leq k\leq k_n}\left\lVert X_{nk}\right\rVert^2:n\in\NN\right\}$,
\begin{displaymath}
\sum_{k=1}^{k_n}E\left(\left\lVert X_{nk}\right\rVert^21_{\left\{\left\lVert X_{nk}\right\rVert\geq\varepsilon\right\}}|\F_{n,k-1}\right)
\rightarrow0\quad\text{in probability as }n\to\infty\text{ for all }\varepsilon>0
\end{displaymath}
follows from Lemma 6.13 in \cite{HaeuslerLuschgy}.
For the proof of (\ref{cond:M1})\,$\Rightarrow$\,(\ref{cond:CLB1}) replace $\left\lVert X_{nk}\right\rVert^2$ 
everywhere by $\left\lVert X_{nk}\right\rVert$ 
and condition (\ref{cond:M2}) by (\ref{cond:M1}).\hfill$\Box$
\medskip

\textit{Proof of Lemma \ref{Lemma}.} Set $U_{nk}=X_{nk,1}$ and $V_{nk}=X_{nk,2}$. 
Let $n\in\NN$ and $\varepsilon>0$ be arbitrary. Then
\begin{eqnarray*}
&&\left|\sum_{k=1}^{k_n}E\left(U_{nk}V_{nk}|\F_{n,k-1}\right)-\sum_{k=1}^{k_n}U_{nk}V_{nk}\right|\\
&&\quad=\left|\sum_{k=1}^{k_n}E\left(U_{nk}V_{nk}1_{\{\left\lVert X_{nk}\right\rVert>\varepsilon\}}|\F_{n,k-1}\right)+
\sum_{k=1}^{k_n}E\left(U_{nk}V_{nk}1_{\{\left\lVert X_{nk}\right\rVert\leq\varepsilon\}}|\F_{n,k-1}\right)\right.\\
&&\quad\left.-\sum_{k=1}^{k_n}U_{nk}V_{nk}1_{\{\left\lVert X_{nk}\right\rVert\leq\varepsilon\}}
-\sum_{k=1}^{k_n}U_{nk}V_{nk}1_{\{\left\lVert X_{nk}\right\rVert>\varepsilon\}}\right|\\
&&\quad\leq\sum_{k=1}^{k_n}E\left(\left|U_{nk}V_{nk}\right|
1_{\{\left\lVert X_{nk}\right\rVert>\varepsilon\}}|\F_{n,k-1}\right)+
\sum_{k=1}^{k_n}\left|U_{nk}V_{nk}\right|1_{\{\left\lVert X_{nk}\right\rVert>\varepsilon\}}\\
&&\quad+\left|\sum_{k=1}^{k_n}U_{nk}V_{nk}1_{\{\left\lVert X_{nk}\right\rVert\leq\varepsilon\}}-
E\left(U_{nk}V_{nk}1_{\{\left\lVert X_{nk}\right\rVert\leq\varepsilon\}}|\F_{n,k-1}\right)\right|\\
&&\quad=:I_n\left(\varepsilon\right)+\II_n\left(\varepsilon\right)+\III_n\left(\varepsilon\right)\,,
\end{eqnarray*}
Note that
\begin{align*}
E\left(I_n\left(\varepsilon\right)\right)&=E\left(\II_n\left(\varepsilon\right)\right)=
\sum_{k=1}^{k_n}E\left(\left|U_{nk}\right|\left|V_{nk}\right|1_{\{\left\lVert X_{nk}\right\rVert>\varepsilon\}}\right)\\
&\leq\sum_{k=1}^{k_n}E\left(\left\lVert X_{nk}\right\rVert^21_{\{\left\lVert X_{nk}\right\rVert>\varepsilon\}}\right)\,.
\end{align*}
Therefore, (\ref{eq:CLBL2Znk}) implies $E\left(I_n\left(\varepsilon\right)\right)\rightarrow0$ and 
$E\left(\II_n\left(\varepsilon\right)\right)\rightarrow0$
as $n\to\infty$ for all $\varepsilon>0$. Moreover, for all $n\in\NN$ and $\varepsilon>0$, because the real-valued
random variables of the bounded martingale difference sequence 
\begin{displaymath}
\left(U_{nk}V_{nk}1_{\{\left\lVert X_{nk}\right\rVert\leq\varepsilon\}}-
E\left(U_{nk}V_{nk}1_{\{\left\lVert X_{nk}\right\rVert\leq\varepsilon\}}|\F_{n,k-1}\right)\right)_{1\leq k\leq k_n}
\end{displaymath}
are pairwise uncorrelated,
\begin{align*}
E\left(\III_n\left(\varepsilon\right)^2\right)&=\sum_{k=1}^{k_n}
E\left[\left(U_{nk}V_{nk}1_{\{\left\lVert X_{nk}\right\rVert\leq\varepsilon\}}-
E\left(U_{nk}V_{nk}1_{\{\left\lVert X_{nk}\right\rVert\leq\varepsilon\}}|\F_{n,k-1}\right)\right)^2\right]\\
&\leq\sum_{k=1}^{k_n}E\left(U_{nk}^2V_{nk}^21_{\{\left\lVert X_{nk}\right\rVert\leq\varepsilon\}}\right)
\leq\varepsilon^2\sum_{k=1}^{k_n}E\left(\left|U_{nk}\right|\left|V_{nk}\right|\right)\\
&\leq\varepsilon^2\sum_{k=1}^{k_n}E\left(\left\lVert X_{nk}\right\rVert^2\right)\leq\varepsilon^2C
\end{align*}
by condition (\ref{eq:L2boundZnk}). Hence for all $n\in\NN$ and $\varepsilon>0$
\begin{displaymath}
E\left(\left|\sum_{k=1}^{k_n}E\left(U_{nk}V_{nk}|\F_{n,k-1}\right)-\sum_{k=1}^{k_n}U_{nk}V_{nk}\right|\right)
\leq E\left(I_n\left(\varepsilon\right)\right)+E\left(\II_n\left(\varepsilon\right)\right)+
\varepsilon C^{1/2}\,.
\end{displaymath}
Because $\varepsilon>0$ is arbitrary and $E\left(I_n\left(\varepsilon\right)\right)\rightarrow0$ and 
$E\left(\II_n\left(\varepsilon\right)\right)\rightarrow0$ as $n\to\infty$ are true, 
the Lemma is proven.\hfill$\Box$
\medskip

\textit{Proof of Proposition \ref{Proposition:CLBandNimplyR}.} Observe that because of
\begin{displaymath}
E\left(\left\lVert X_{nk}\right\rVert
1_{\left\{\left\lVert X_{nk}\right\rVert\geq\varepsilon\right\}}|\F_{n,k-1}\right)\leq\frac{1}{\varepsilon}
E\left(\left\lVert X_{nk}\right\rVert^21_{\left\{\left\lVert X_{nk}\right\rVert\geq\varepsilon\right\}}|\F_{n,k-1}\right)
\end{displaymath}
for square integrable arrays condition (\ref{cond:CLB1}) is an immediate consequence of condition (\ref{cond:CLB}).
\smallskip

Condition \eqref{cond:R} follows from \eqref{cond:N} and
\begin{displaymath}
\sum_{k=1}^{k_n}E\left(X_{nk}X_{nk}^T|\F_{n,k-1}\right)-\sum_{k=1}^{k_n}X_{nk}X_{nk}^T\rightarrow
0\quad\text{in probability as }n\to\infty
\end{displaymath}
which is equivalent to
\begin{equation}
\label{eq:differencePropCLBandNimplyR}
\sum_{k=1}^{k_n}E\left(X_{nk,j}X_{nk,\ell}|\F_{n,k-1}\right)-\sum_{k=1}^{k_n}X_{nk,j}X_{nk,\ell}
\rightarrow0\quad\text{in probability as }n\to\infty
\end{equation}
for all $1\leq j,\ell\leq d$. For the proof of \eqref{eq:differencePropCLBandNimplyR} 
let $j,\ell$ be arbitrary. Set $U_{nk}=X_{nk,j}$ and
$V_{nk}=X_{nk,\ell}$. Then we have to show
\begin{equation}
\label{eq:UVdifferencePropCLBandNimplyR}
\sum_{k=1}^{k_n}E\left(U_{nk}V_{nk}|\F_{n,k-1}\right)-\sum_{k=1}^{k_n}U_{nk}V_{nk}\rightarrow
0\quad\text{in probability as }n\to\infty\,.
\end{equation}
For the proof of (\ref{eq:UVdifferencePropCLBandNimplyR}) for each $0<c<\infty$ and all 
$n\in\NN$ we define the stopping times
\begin{align*}
&\tau_n^U\left(c\right)=\max\left\{k\in\{0,1,\ldots,k_n\}:\sum_{i=1}^kE\left(U_{ni}^2|\F_{n,i-1}\right)
\leq c\right\}\,,\\
&\tau_n^V\left(c\right)=\max\left\{k\in\{0,1,\ldots,k_n\}:\sum_{i=1}^kE\left(V_{ni}^2|\F_{n,i-1}\right)
\leq c\right\}\qquad\text{and}\\
&\tau_n\left(c\right)=\tau_n^U\left(c\right)\wedge\tau_n^V\left(c\right)
\end{align*}
w.r.t. $\left(\F_{nk}\right)_{0\leq k\leq k_n}$ with $\sum_\emptyset=0$ and the real-valued random variables
\begin{displaymath}
U_{nk}\left(c\right)=U_{nk}1_{\{k\leq\tau_n\left(c\right)\}}\quad\text{and}\quad
V_{nk}\left(c\right)=V_{nk}1_{\{k\leq\tau_n\left(c\right)\}}
\end{displaymath}
for all $n\in\NN$ and $1\leq k\leq k_n$. Observe that these random variables are square integrable because 
$\left|U_{nk}\left(c\right)\right|\leq\left|U_{nk}\right|$ 
and $\left|V_{nk}\left(c\right)\right|\leq\left|V_{nk}\right|$.
Then for all $n\in\NN$ and $0<c<\infty$
\begin{eqnarray}
&&\sum_{k=1}^{k_n}E\left(U_{nk}V_{nk}|\F_{n,k-1}\right)-\sum_{k=1}^{k_n}U_{nk}V_{nk}
\label{eq:UVdifferencerepPropCLBandNimplyR}\\
&&\quad=\sum_{k=1}^{k_n}E\left(U_{nk}V_{nk}|\F_{n,k-1}\right)-\sum_{k=1}^{k_n}E\left(U_{nk}\left(c\right)V_{nk}
\left(c\right)|\F_{n,k-1}\right)\nonumber\\
&&\quad\quad+\sum_{k=1}^{k_n}E\left(U_{nk}\left(c\right)V_{nk}\left(c\right)|\F_{n,k-1}\right)-
\sum_{k=1}^{k_n}U_{nk}\left(c\right)V_{nk}\left(c\right)\nonumber\\
&&\quad\quad+\sum_{k=1}^{k_n}U_{nk}\left(c\right)V_{nk}\left(c\right)-\sum_{k=1}^{k_n}U_{nk}V_{nk}\nonumber\\
&&\quad=:I_n\left(c\right)+\II_n\left(c\right)+\III_n\left(c\right)\,.\nonumber
\end{eqnarray}
Now, using that $\tau_n\left(c\right)$ is a stopping time w.r.t. 
$\left(\F_{nk}\right)_{1\leq k\leq k_n}$, we have
\begin{align*}
I_n\left(c\right)&=\sum_{k=1}^{k_n}E\left(U_{nk}V_{nk}|\F_{n,k-1}\right)-
\sum_{k=1}^{\tau_n\left(c\right)}E\left(U_{nk}V_{nk}|\F_{n,k-1}\right)\\
&=\sum_{k=\tau_n\left(c\right)+1}^{k_n}E\left(U_{nk}V_{nk}|\F_{n,k-1}\right)
\end{align*}
so that for all $\varepsilon>0$
\begin{eqnarray}
&&P\left(\left|I_n\left(c\right)\right|\geq\varepsilon\right)\leq P\left(\tau_n\left(c\right)<k_n\right)
\leq P\left(\tau_n^U\left(c\right)<k_n\right)+P\left(\tau_n^V\left(c\right)<k_n\right)
\label{eq:boundInPropCLBandNimplyR}\\
&&\quad\leq P\left(\sum_{i=1}^{k_n}E\left(U_{ni}^2|\F_{n,i-1}\right)>c\right)+
P\left(\sum_{i=1}^{k_n}E\left(V_{ni}^2|\F_{n,i-1}\right)>c\right)\nonumber\,.
\end{eqnarray}
Similarly
\begin{eqnarray}
&&P\left(\left|\III_n\left(c\right)\right|\geq\varepsilon\right)\label{eq:boundIIInPropCLBandNimplyR}\\
&&\quad\leq P\left(\sum_{i=1}^{k_n}E\left(U_{ni}^2|\F_{n,i-1}\right)>c\right)+
P\left(\sum_{i=1}^{k_n}E\left(V_{ni}^2|\F_{n,i-1}\right)>c\right)\,.\nonumber
\end{eqnarray}
Next we will show $\II_n\left(c\right)\rightarrow0$ in $L_1$ as $n\to\infty$ for all $0<c<\infty$ by an
application of Lemma \ref{Lemma}. For this, we set
$Z_{nk}\left(c\right)=\left(U_{nk}\left(c\right),V_{nk}\left(c\right)\right)^T$
for all $n\in\NN$ and $1\leq k\leq k_n$. Note that for $0<c<\infty$ and $n\in\NN$
\begin{eqnarray}
&&\sum_{k=1}^{k_n}E\left(\left\lVert Z_{nk}\left(c\right)\right\rVert^2|\F_{n,k-1}\right)=
\sum_{k=1}^{k_n}E\left(U_{kn}^2\left(c\right)+V_{nk}^2\left(c\right)|\F_{n,k-1}\right)
\label{eq:boundZnkPropCLBandNimplyR}\\
&&\quad=\sum_{k=1}^{\tau_n\left(c\right)}E\left(U_{kn}^2|\F_{n,k-1}\right)+
\sum_{k=1}^{\tau_n\left(c\right)}E\left(V_{nk}^2|\F_{n,k-1}\right)\nonumber\\
&&\quad\leq\sum_{k=1}^{\tau_n^U\left(c\right)}E\left(U_{kn}^2|\F_{n,k-1}\right)+
\sum_{k=1}^{\tau_n^V\left(c\right)}E\left(V_{nk}^2|\F_{,k-1}\right)\leq 2c\nonumber
\end{eqnarray}
by definiton of $\tau_n^U$ and $\tau_n^V$. This implies
\begin{displaymath}
\sum_{k=1}^{k_n}E\left(\left\lVert Z_{nk}\left(c\right)\right\rVert^2\right)\leq2c
\end{displaymath}
so that condition \eqref{eq:L2boundZnk} is satisfied, Observe also that
\begin{displaymath}
\left\lVert Z_{nk}\left(c\right)\right\rVert^2=U_{kn}^2\left(c\right)+V_{nk}^2\left(c\right)\leq U_{kn}^2+V_{nk}^2=
X_{kn,j}^2+X_{nk,\ell}^2\leq2\left\lVert X_{nk}\right\rVert^2
\end{displaymath}
where the factor $2$ is needed only in case $j=\ell$ so that
\begin{displaymath}
\sum_{k=1}^{k_n}E\left(\left\lVert Z_{nk}\left(c\right)\right\rVert^2
1_{\{\left\lVert Z_{nk}\left(c\right)\right\rVert\geq\varepsilon\}}|\F_{n,k-1}\right)
\leq2\sum_{k=1}^{k_n}E\left(\left\lVert X_{nk}\right\rVert^2
1_{\{\left\lVert X_{nk}\right\rVert\geq\varepsilon/2^{1/2}\}}|\F_{n,k-1}\right)
\end{displaymath}
which by \eqref{cond:CLB} implies
\begin{displaymath}
\sum_{k=1}^{k_n}E\left(\left\lVert Z_{nk}\left(c\right)\right\rVert^2
1_{\{\left\lVert Z_{nk}\left(c\right)\right\rVert\geq\varepsilon\}}|\F_{n,k-1}\right)
\rightarrow0\quad\text{in probability as }n\to\infty
\end{displaymath}
for every $\varepsilon>0$. From this and the bound in \eqref{eq:boundZnkPropCLBandNimplyR} we obtain
\begin{displaymath}
\sum_{k=1}^{k_n}E\left(\left\lVert Z_{nk}\left(c\right)\right\rVert^2
1_{\{\left\lVert Z_{nk}\left(c\right)\right\rVert\geq\varepsilon\}}\right)
\rightarrow0\quad\text{as }n\to\infty\text{ for every }\varepsilon>0\
\end{displaymath}
by dominated convergence. Thus, \eqref{eq:CLBL2Znk} is also satisfied and 
Lemma~\ref{Lemma} yields $\II_n\left(c\right)\rightarrow0$
in $L_1$ as $n\to\infty$.  From (\ref{eq:UVdifferencerepPropCLBandNimplyR}) 
we get for all $\varepsilon>0$ and $0<c<\infty$, 
using \eqref{eq:boundInPropCLBandNimplyR} 
and \eqref{eq:boundIIInPropCLBandNimplyR},
\begin{eqnarray*}
&&P\left(\left|\sum_{k=1}^{k_n}E\left(U_{nk}V_{nk}|\F_{n,k-1}\right)-\sum_{k=1}^{k_n}U_{nk}V_{nk}
\right|\geq3\varepsilon\right)\\
&&\quad\leq P\left(\left|I_n\left(c\right)\right|\geq\varepsilon\right)+
P\left(\left|\II_n\left(c\right)\right|\geq\varepsilon\right)+
P\left(\left|\III_n\left(c\right)\right|\geq\varepsilon\right)\\
&&\quad\leq 2P\left(\sum_{i=1}^{k_n}E\left(U_{ni}^2|\F_{n,i-1}\right)>c\right)+
2P\left(\sum_{i=1}^{k_n}E\left(V_{ni}^2|\F_{n,i-1}\right)>c\right)\\
&&\qquad\qquad+P\left(\left|\II_n\left(c\right)\right|\geq\varepsilon\right)\\
&&\quad=2P\left(\sum_{i=1}^{k_n}E\left(X_{ni,j}^2|\F_{n,i-1}\right)>c\right)+
2P\left(\sum_{i=1}^{k_n}E\left(X_{ni,\ell}^2|\F_{n,i-1}\right)>c\right)\\
&&\qquad\qquad+P\left(\left|\II_n\left(c\right)\right|\geq\varepsilon\right)\,.
\end{eqnarray*}
Now $P\left(\left|\II_n\left(c\right)\right|\geq\varepsilon\right)\rightarrow0$ 
as $n\to\infty$ for each $\varepsilon>0$ and $0<c<\infty$ because of
$\II_n\left(c\right)\rightarrow0$ in $L_1$, and the sequences
 $\left(\sum_{i=1}^{k_n}E\left(X_{ni,j}^2|\F_{n,i-1}\right)\right)_{n\in\NN}$ and
$\left(\sum_{i=1}^{k_n}E\left(X_{ni,\ell}^2|\F_{n,i-1}\right)\right)_{n\in\NN}$ are 
bounded in probability as a consequence of condition
(\ref{cond:N}), from which the assertion follows.\hfill$\Box$
\medskip

\textit{Proof of Proposition~\ref{Proposition:CLBandRimplyN}.}
Condition (\ref{cond:N}) follows from (\ref{cond:R}) and
\begin{displaymath}
\sum_{k=1}^{k_n}E\left(X_{nk}X_{nk}^T|\F_{n,k-1}\right)-\sum_{k=1}^{k_n}X_{nk}X_{nk}^T\rightarrow0
\quad\text{in probability as }n\to\infty
\end{displaymath}
which is equivalent to
\begin{equation}
\label{eq:difference}
\sum_{k=1}^{k_n}E\left(X_{nk,j}X_{nk,\ell}|\F_{n,k-1}\right)-\sum_{k=1}^{k_n}X_{nk,j}X_{nk,\ell}\rightarrow0
\quad\text{in probability as }n\to\infty
\end{equation}
for all $1\leq j,\ell\leq d$. For the proof of (\ref{eq:difference}) let $j$ and $\ell$ 
be fixed and set $U_{nk}=X_{nk,j}$
and $V_{nk}=X_{nk,\ell}$. Then we have to show
\begin{equation}
\label{eq:UVdifferencePropApprox}
\sum_{k=1}^{k_n}E\left(U_{nk}V_{nk}|\F_{n,k-1}\right)-
\sum_{k=1}^{k_n}U_{nk}V_{nk}\rightarrow0\quad\text{in probability as }n\to\infty\,.
\end{equation}
For the proof of (\ref{eq:UVdifferencePropApprox}) for each $0<c<\infty$ and all $n\in\NN$ 
we define the stopping times
\begin{align*}
&\tau_n^U\left(c\right)=\min\left\{k\in\{1,\ldots,k_n\}:\sum_{i=1}^kU_{ni}^2>c\right\}\wedge k_n\,,\\
&\tau_n^V\left(c\right)=\min\left\{k\in\{1,\ldots,k_n\}:\sum_{i=1}^kV_{ni}^2>c\right\}\wedge k_n\qquad\text{and}\\
&\tau_n\left(c\right)=\tau_n^U\left(c\right)\wedge\tau_n^V\left(c\right)
\end{align*}
w.r.t. $\left(\F_{nk}\right)_{0\leq k\leq k_n}$ with $\min\emptyset=\infty$ and the real-valued random variables
\begin{displaymath}
U_{nk}\left(c\right)=U_{nk}1_{\{k\leq\tau_n\left(c\right)\}}\quad\text{and}\quad
V_{nk}\left(c\right)=V_{nk}1_{\{k\leq\tau_n\left(c\right)\}}
\end{displaymath}
for all $n\in\NN$ and $1\leq k\leq k_n$. Then for all $n\in\NN$
\begin{eqnarray}
&&\sum_{k=1}^{k_n}E\left(U_{nk}V_{nk}|\F_{n,k-1}\right)-\sum_{k=1}^{k_n}U_{nk}V_{nk}
\label{eq:UVdifferencerepPropApprox}\\
&&\quad=\sum_{k=1}^{k_n}E\left(U_{nk}V_{nk}|\F_{n,k-1}\right)-\sum_{k=1}^{k_n}E\left(U_{nk}\left(c\right)V_{nk}
\left(c\right)|\F_{n,k-1}\right)\nonumber\\
&&\quad\quad+\sum_{k=1}^{k_n}E\left(U_{nk}\left(c\right)V_{nk}\left(c\right)|\F_{n,k-1}\right)-
\sum_{k=1}^{k_n}U_{nk}\left(c\right)V_{nk}\left(c\right)\nonumber\\
&&\quad\quad+\sum_{k=1}^{k_n}U_{nk}\left(c\right)V_{nk}\left(c\right)-\sum_{k=1}^{k_n}U_{nk}V_{nk}\nonumber\\
&&\quad=:I_n\left(c\right)+\II_n\left(c\right)+\III_n\left(c\right)\,.\nonumber
\end{eqnarray}
Now, using that $\tau_n\left(c\right)$ is a stopping time w.r.t. 
$\left(\F_{nk}\right)_{0\leq k\leq k_n}$, we have
\begin{align*}
I_n\left(c\right)&=\sum_{k=1}^{k_n}E\left(U_{nk}V_{nk}|\F_{n,k-1}\right)-
\sum_{k=1}^{\tau_n\left(c\right)}E\left(U_{nk}V_{nk}|\F_{n,k-1}\right)\\
&=\sum_{k=\tau_n\left(c\right)+1}^{k_n}E\left(U_{nk}V_{nk}|\F_{n,k-1}\right)
\end{align*}
so that for all $\varepsilon>0$
\begin{eqnarray}
&&P\left(\left|I_n\left(c\right)\right|\geq\varepsilon\right)\leq P\left(\tau_n\left(c\right)<k_n\right)
\leq P\left(\tau_n^U\left(c\right)<k_n\right)+P\left(\tau_n^V\left(c\right)<k_n\right)
\label{eq:boundInPropApprox}\\
&&\quad\leq P\left(\sum_{i=1}^{k_n}U_{ni}^2>c\right)+P\left(\sum_{i=1}^{k_n}V_{ni}^2>c\right)\nonumber\,.
\end{eqnarray}
Similarly
\begin{equation}\label{eq:boundIIInPropApprox}
P\left(\left|\III_n\left(c\right)\right|\geq\varepsilon\right)\leq
P\left(\sum_{i=1}^{k_n}U_{ni}^2>c\right)+P\left(\sum_{i=1}^{k_n}V_{ni}^2>c\right)\,.
\end{equation}
Next we will show $\II_n\left(c\right)\rightarrow0$ in $L_1$ as $n\to\infty$ for all $0<c<\infty$ by an
application of Lemma \ref{Lemma}. For this, we set
 $Z_{nk}\left(c\right)=\left(U_{nk}\left(c\right),V_{nk}\left(c\right)\right)^T$
for all $n\in\NN$ and $1\leq k\leq k_n$. Note that for $0<c<\infty$ and $n\in\NN$
\begin{eqnarray}
&&\sum_{k=1}^{k_n}\left\lVert 
Z_{nk}\left(c\right)\right\rVert^2=\sum_{k=1}^{k_n}U_{kn}^2\left(c\right)+V_{nk}^2\left(c\right)=
\sum_{k=1}^{\tau_n\left(c\right)}U_{kn}^2+\sum_{k=1}^{\tau_n\left(c\right)}V_{nk}^2\label{eq:boundZnkPropApprox}\\
&&\quad\leq\sum_{k=1}^{\tau_n^U\left(c\right)}U_{kn}^2+\sum_{k=1}^{\tau_n^V\left(c\right)}V_{nk}^2\nonumber\\
&&\quad\leq\sum_{k=1}^{\tau_n^U\left(c\right)-1}U_{kn}^2+\max_{1\leq k\leq k_n}U_{nk}^2+
\sum_{k=1}^{\tau_n^V\left(c\right)-1}V_{nk}^2+\max_{1\leq k\leq k_n}V_{nk}^2\leq2c+2a^2\nonumber
\end{eqnarray}
by definiton of $\tau_n^U$ and $\tau_n^V$ and because $U_{nk}^2\leq\left\lVert X_{nk}\right\rVert^2\leq a^2$ and likewise
$V_{nk}^2\leq a^2$ for all $n\in\NN$ and $1\leq k\leq k_n$. Consequently,
\begin{displaymath}
\sum_{k=1}^{k_n}E\left(\left\lVert Z_{nk}\left(c\right)\right\rVert^2\right)=
E\left(\sum_{k=1}^{k_n}\left\lVert Z_{nk}\left(c\right)\right\rVert^2\right)\leq2c+2a^2\,,
\end{displaymath}
which shows that $\left(Z_{nk}\left(c\right)\right)_{1\leq k\leq k_n,n\in\NN}$ 
satisfies condition (\ref{eq:L2boundZnk})
in Lemma~\ref{Lemma}. Next, we will show that condition (\ref{eq:CLBL2Znk}) is also satisfied. 
For this, note that for $n\in\NN$
and $1\leq k\leq k_n$
\begin{displaymath}
\left\lVert Z_{nk}\left(c\right)\right\rVert^2=U_{nk}^2\left(c\right)+V_{nk}^2\left(c\right)
\leq U_{nk}^2+V_{nk}^2=X_{nk,j}^2+
X_{nk,\ell}^2\leq2\left\lVert X_{nk}\right\rVert^2 
\end{displaymath}
where the factor $2$ is needed only in case $j=\ell$ so that
\begin{displaymath}
\sum_{k=1}^{k_n}E\left(\left\lVert Z_{nk}\left(c\right)\right\rVert^2
1_{\{\left\lVert Z_{nk}\left(c\right)\right\rVert\geq\varepsilon\}}|\F_{n,k-1}\right)
\leq2\sum_{k=1}^{k_n}E\left(\left\lVert X_{nk}\right\rVert^2
1_{\{\left\lVert X_{nk}\right\rVert\geq\varepsilon/2^{1/2}\}}|\F_{n,k-1}\right)\,.
\end{displaymath}
Therefore, (\ref{cond:CLB}) and Lemma 6.5 in \cite{HaeuslerLuschgy} imply
\begin{displaymath}
\sum_{k=1}^{k_n}\left\lVert Z_{nk}\left(c\right)\right\rVert^2
1_{\{\left\lVert Z_{nk}\left(c\right)\right\rVert\geq\varepsilon\}}\rightarrow0
\quad\text{in probability as }n\to\infty\text{ for every }\varepsilon>0\,.
\end{displaymath}
In view of (\ref{eq:boundZnkPropApprox}) condition (\ref{eq:CLBL2Znk}) 
follows by dominated convergence. Consequently, from Lemma \ref{Lemma} we obtain
$\II_n\left(c\right)\rightarrow0$ in $L_1$ as $n\to\infty$
for all $0<c<\infty$. From (\ref{eq:UVdifferencerepPropApprox}) we get for all 
$\varepsilon>0$ and $0<c<\infty$, using \eqref{eq:boundInPropApprox} 
and \eqref{eq:boundIIInPropApprox},
\begin{eqnarray*}
&&P\left(\left\lvert\sum_{k=1}^{k_n}E\left(U_{nk}V_{nk}\rvert\F_{n,k-1}\right)-\sum_{k=1}^{k_n}U_{nk}V_{nk})
\right|\geq3\varepsilon\right)\\
&&\quad\leq 
P\left(\left\lvert I_n\left(c\right)\right\rvert\geq\varepsilon\right)+P\left(\left\lvert\II_n\left(c\right)
\right\rvert\geq\varepsilon\right)+
P\left(\left\lvert\III_n\left(c\right)\right\rvert\geq\varepsilon\right)\\
&&\quad\leq 2P\left(\sum_{i=1}^{k_n}U_{ni}^2>c\right)+2P\left(\sum_{i=1}^{k_n}V_{ni}^2>c\right)+
P\left(\left\lvert\II_n\left(c\right)\right\vert\geq\varepsilon\right)\\
&&\quad=2P\left(\sum_{i=1}^{k_n}X_{ni,j}^2>c\right)+2P\left(\sum_{i=1}^{k_n}X_{ni,\ell}^2>c\right)+
P\left(\left\lvert\II_n\left(c\right)\right\rvert\geq\varepsilon\right)\,.
\end{eqnarray*}
Now $P\left(\left\lvert\II_n\left(c\right)\right\rvert\geq\varepsilon\right)\rightarrow0$ as $n\to\infty$ 
for each $\varepsilon>0$ and $0<c<\infty$ because of
$\II_n\left(c\right)\rightarrow0$ in $L_1$, and the sequences 
$\left(\sum_{i=1}^{k_n}X_{ni,j}^2\right)_{n\in\NN}$ and
$\left(\sum_{i=1}^{k_n}X_{ni,\ell}^2\right)_{n\in\NN}$ are bounded in probability as a consequence of condition
(\ref{cond:R}), from which the assertion follows.\hfill$\Box$
\medskip

\textit{Proof of Proposition \ref{Proposition:CLB1andRimplytheT}.} 
Let $a>0$ be arbitrary. From the martingale difference property of 
$\left(X_{nk}\right)_{1\leq k\leq k_n,n\in\NN}$ we obtain
for all $n\in\NN$ and $1\leq k\leq k_n$
\begin{equation}
\label{eq:martingaleequation}
E\left(X_{nk}1_{\left\{\left\lVert X_{nk}\right\rVert\leq a\right\}}|\F_{n,k-1}\right)
=-E\left(X_{nk}1_{\left\{\left\lVert X_{nk}\right\rVert>a\right\}}|\F_{n,k-1}\right)\,.
\end{equation}
Proof of (\ref{cond:Ta}): By (\ref{eq:martingaleequation})
\begin{eqnarray*}
&&\left\lVert\sum_{k=1}^{k_n}E\left(X_{nk}1_{\{\left\lVert X_{nk}\right\rVert\leq a\}}|\F_{n,k-1}\right)\right\rVert\leq
\sum_{k=1}^{k_n}\left\lVert E\left(X_{nk}1_{\{\left\lVert X_{nk}\right\rVert
\leq a\}}|\F_{n,k-1}\right)\right\rVert\\
&&\quad=\sum_{k=1}^{k_n}\left\lVert E\left(X_{nk}
1_{\{\left\lVert X_{nk}\right\rVert>a\}}|\F_{n,k-1}\right)\right\rVert\leq
\sum_{k=1}^{k_n}E\left(\left\lVert X_{nk}\right\rVert1_{\{\left\lVert X_{nk}\right\rVert>a\}}|\F_{n,k-1}\right)\,,
\end{eqnarray*}
and from (\ref{cond:CLB1}) we obtain
\begin{equation}
\label{eq:secondhalf}
\sum_{k=1}^{k_n}E\left(X_{nk}1_{\{\left\lVert X_{nk}\right\rVert\leq a\}}|\F_{n,k-1}\right)\rightarrow
0\quad\text{in probability as }n\to\infty\,.
\end{equation}
By Lemma 6.5 in \cite{HaeuslerLuschgy} it follows from \eqref{cond:CLB1} that
\begin{displaymath}
\sum_{k=1}^{k_n}\left\lVert X_{nk}\right\rVert1_{\{\left\lVert X_{nk}\right\rVert>a\}}\rightarrow
0\quad\text{in probability as }n\to\infty\,
\end{displaymath}
and therefore also
\begin{equation}
\label{eq:firsthalf}
\sum_{k=1}^{k_n}X_{nk}1_{\{\left\lVert X_{nk}\right\rVert>a\}}\rightarrow0\quad\text{in probability as }n\to\infty\,.
\end{equation}
Combining (\ref{eq:firsthalf}) and (\ref{eq:secondhalf}) we get (\ref{cond:Ta}).
\smallskip

Proof of (\ref{cond:TMa}): For all $n\in\NN$ by (\ref{eq:martingaleequation})
\begin{align*}
\max_{1\leq k\leq k_n}\left\lVert X_{nk}\left(a\right)\right\rVert
&\leq\max_{1\leq k\leq k_n}\left\lVert X_{nk}\right\rVert+
\max_{1\leq k\leq k_n}\left\lVert E\left(X_{nk}1_{\{\left\lVert X_{nk}\right\rVert
\leq a\}}|\F_{n,k-1}\right)\right\rVert\\
&=\max_{1\leq k\leq k_n}\left\lVert X_{nk}\right\rVert+
\max_{1\leq k\leq k_n}\left\lVert E\left(X_{nk}
1_{\{\left\lVert X_{nk}\right\rVert>a\}}|\F_{n,k-1}\right)\right\rVert\\
&\leq\max_{1\leq k\leq k_n}\left\lVert X_{nk}\right\rVert+
\sum_{k=1}^{k_n}E\left(\left\lVert X_{nk}\right\rVert1_{\{\left\lVert X_{nk}\right\rVert>a\}}|\F_{n,k-1}\right)\,.
\end{align*}
Note that the second summand on the right-hand side converges to zero in probability as 
$n\to\infty$ by (\ref{cond:CLB1}) and that this condition
also implies 
\begin{equation}
\label{eq:convergencemax}
\max_{1\leq k\leq k_n}\left\lVert X_{nk}\right\rVert\rightarrow0\quad\text{in probability as }n\to\infty
\end{equation}
in view of Proposition 6.6 in \cite{HaeuslerLuschgy} which completes the proof.
\smallskip

Proof of (\ref{cond:TRa}): Condition (\ref{cond:TRa}) follows from (\ref{cond:R}) and
\begin{displaymath}
\sum_{k=1}^{k_n}X_{nk}X_{nk}^T-\sum_{k=1}^{k_n}X_{nk}\left(a\right)X_{nk}\left(a\right)^T\rightarrow
0\quad\text{in probability as }n\to\infty\,,
\end{displaymath}
which is equivalent to
\begin{equation}
\label{eq:differenceXnkXnka}
\sum_{k=1}^{k_n}X_{nk,j}X_{nk,\ell}-\sum_{k=1}^{k_n}X_{nk,j}\left(a\right)X_{nk,\ell}\left(a\right)
\rightarrow0\quad\text{in probability as }n\to\infty
\end{equation}
for all $1\leq j,\ell\leq d$. For the proof of (\ref{eq:differenceXnkXnka}) let $j,\ell$ be arbitrary. 
Then for all $n\in\NN$
\begin{eqnarray*}
&&\sum_{k=1}^{k_n}X_{nk,j}X_{nk,\ell}-\sum_{k=1}^{k_n}X_{nk,j}\left(a\right)X_{nk,\ell}\left(a\right)\\
&&\quad=\sum_{k=1}^{k_n}X_{nk,j}X_{nk,\ell}-\sum_{k=1}^{k_n}\left[X_{nk,j}
1_{\{\left\Vert X_{nk}\right\Vert\leq a\}}-
E\left(X_{nk,j}1_{\{\left\Vert X_{nk}\right\Vert\leq a\}}|\F_{n,k-1}\right)\right]\times\\
&&\hspace{126pt}\left[X_{nk,\ell}1_{\{\left\Vert X_{nk}\right\Vert\leq a\}}-E\left(X_{nk,\ell}
1_{\{\left\Vert X_{nk}\right\Vert\leq a\}}|\F_{n,k-1}\right)\right]\\
&&\quad=\sum_{k=1}^{k_n}X_{nk,j}X_{nk,\ell}-\sum_{k=1}^{k_n}X_{nk,j}X_{nk,\ell}
1_{\{\left\Vert X_{nk}\right\Vert\leq a\}}\\
&&\quad\quad+\sum_{k=1}^{k_n}X_{nk,j}1_{\{\left\Vert X_{nk}\right\Vert\leq a\}}E\left(X_{nk,\ell}
1_{\{\left\Vert X_{nk}\right\Vert\leq a\}}|\F_{n,k-1}\right)\\
&&\quad\quad+\sum_{k=1}^{k_n}X_{nk,\ell}1_{\{\left\Vert X_{nk}\right\Vert\leq a\}}E\left(X_{nk,j}
1_{\{\left\Vert X_{nk}\right\Vert\leq a\}}|\F_{n,k-1}\right)\\
&&\quad\quad-\sum_{k=1}^{k_n}E\left(X_{nk,j}1_{\{\left\Vert X_{nk}\right\Vert\leq a\}}|\F_{n,k-1}\right)
E\left(X_{nk,\ell}1_{\{\left\Vert X_{nk}\right\Vert\leq a\}}|\F_{n,k-1}\right)\\
&&\quad=:I_n+\II_n+\III_n-\IV_n\,.
\end{eqnarray*}
We have
\begin{align*}
\left|I_n\right|&=\left|\sum_{k=1}^{k_n}X_{nk,j}X_{nk,\ell}1_{\{\left\Vert X_{nk}\right\Vert>a\}}\right|\leq
\sum_{k=1}^{k_n}\left|X_{nk,j}\right|\left|X_{nk,\ell}\right|1_{\{\left\Vert X_{nk}\right\Vert>a\}}\\
&\leq
\sum_{k=1}^{k_n}\left\lVert X_{nk}\right\rVert^21_{\{\left\Vert X_{nk}\right\Vert>a\}}\,.
\end{align*}
Now for all $n\in\NN$ and $\varepsilon>0$
\begin{displaymath}
P\left(\sum_{k=1}^{k_n}\left\lVert X_{nk}\right\rVert^2
1_{\{\left\Vert X_{nk}\right\Vert>a\}}\geq\varepsilon\right)\leq
P\left(\max_{1\leq k\leq k_n}\left\lVert X_{nk}\right\rVert>a\right)\rightarrow0\quad\text{as }n\to\infty
\end{displaymath}
by (\ref{eq:convergencemax}) so that
\begin{displaymath}
\sum_{k=1}^{k_n}\left\lVert X_{nk}\right\rVert^21_{\{\left\Vert X_{nk}\right\Vert>a\}}
\rightarrow0\quad\text{in probability as }n\to\infty
\end{displaymath}
which clearly entails $I_n\rightarrow0$ in probability as $n\to\infty$.
\smallskip

For all $n\in\NN$ we have, again using (\ref{eq:martingaleequation}),
\begin{align*}
\left|\II_n\right|&=\left|\sum_{k=1}^{k_n}X_{nk,j}1_{\{\left\Vert X_{nk}\right\Vert\leq a\}}
E\left(X_{nk,\ell}1_{\{\left\Vert X_{nk}\right\Vert\leq a\}}|\F_{n,k-1}\right)\right|\\
&\leq\sum_{k=1}^{k_n}\left|X_{nk,j}\right|1_{\{\left\Vert X_{nk}\right\Vert\leq a\}}
\left|E\left(X_{nk,\ell}1_{\{\left\Vert X_{nk}\right\Vert\leq a\}}|\F_{n,k-1}\right)\right|\\
&=\sum_{k=1}^{k_n}\left|X_{nk,j}\right|1_{\{\left\Vert X_{nk}\right\Vert\leq a\}}
\left|E\left(X_{nk,\ell}1_{\{\left\Vert X_{nk}\right\Vert>a\}}|\F_{n,k-1}\right)\right|\\
&\leq a\sum_{k=1}^{k_n}E\left(\left|X_{nk,\ell}\right|1_{\{\left\Vert X_{nk}\right\Vert>a\}}|\F_{n,k-1}\right)\\
&\leq a\sum_{k=1}^{k_n}E\left(\left\lVert X_{nk}\right\rVert
1_{\{\left\Vert X_{nk}\right\Vert>a\}}|\F_{n,k-1}\right)\,.
\end{align*}
By (\ref{cond:CLB1}) the right-hand side of these inequalities converges to zero 
in probability as $n\to\infty$ so that
also $\II_n\rightarrow0$ in probability as well as $\III_n\rightarrow0$ in probability simply by interchanging
$j$ and $\ell$.
\smallskip

Finally, for all $n\in\NN$, again using (\ref{eq:martingaleequation}),
\begin{align*}
\left|\IV_n\right|&\leq\sum_{k=1}^{k_n}\left|E\left(X_{nk,j}
1_{\{\left\Vert X_{nk}\right\Vert\leq a\}}|\F_{n,k-1}\right)\right|
\left|E\left(X_{nk,\ell}1_{\{\left\Vert X_{nk}\right\Vert\leq a\}}|\F_{n,k-1}\right)\right|\\
&=\sum_{k=1}^{k_n}\left|E\left(X_{nk,j}1_{\{\left\Vert X_{nk}\right\Vert>a\}}|\F_{n,k-1}\right)\right|
\left|E\left(X_{nk,\ell}1_{\{\left\Vert X_{nk}\right\Vert\leq a\}}|\F_{n,k-1}\right)\right|\\
&\leq\sum_{k=1}^{k_n}E\left(\left|X_{nk,j}\right|1_{\{\left\Vert X_{nk}\right\Vert>a\}}|\F_{n,k-1}\right)
E\left(\left|X_{nk,\ell}\right|1_{\{\left\Vert X_{nk}\right\Vert\leq a\}}|\F_{n,k-1}\right)\\
&\leq a\sum_{k=1}^{k_n}E\left(\left\lVert X_{nk}\right\rVert
1_{\{\left\Vert X_{nk}\right\Vert>a\}}|\F_{n,k-1}\right)\,,
\end{align*}
where the right-hand side converges to zero in probability as $n\to\infty$ by 
(\ref{cond:CLB1}) which proves $\IV_n\rightarrow0$ in probability
as $n\to\infty$ and concludes the proof of the proposition.\hfill$\Box$
\medskip

\textit{Proof of Proposition \ref{Proposition:Approx}.} (a) From \eqref{cond:TMa} 
and $\left\lVert X_{nk}\left(a\right)\right\rVert\leq2a$ for all $n$ and $k$ 
we see that \eqref{cond:M1} holds true for the truncated array 
$\left(X_{nk}\left(a\right)\right)_{1\leq k\leq k_n,n\in\NN}$
which clearly is a uniformly bounded martingale difference array w.r.t. 
$\left(\F_{nk}\right)_{0\leq k\leq k_n,n\in\NN}$. Therefore, \eqref{cond:CLB1}
also holds true by Proposition~\ref{Proposition:MimpliesCLB}. 
Assumption~\eqref{cond:TRa} is exactly condition~\eqref{cond:R} for
the truncated array $\left(X_{nk}\left(a\right)\right)_{1\leq k\leq k_n,n\in\NN}$. 
As to~\eqref{eq:approxXnka}, for all $n\in\NN$ we have
\begin{displaymath}
\sum_{k=1}^{k_n}X_{nk}-\sum_{k=1}^{k_n}X_{nk}\left(a\right)\\
=\sum_{k=1}^{k_n}X_{nk}1_{\{\left\lVert X_{nk}\right\rVert>a\}}+\sum_{k=1}^{k_n}E\left(X_{nk}1_{\{\left\lVert X_{nk}\right\rVert
\leq a\}}|\F_{n,k-1}\right)\,,
\end{displaymath}
and (\ref{eq:approxXnka}) is equivalent to (\ref{cond:Ta}).
\smallskip

(b) Condition~\eqref{cond:M2} for the uniformly bounded array 
$\left(X_{nk}\left(a\right)\right)_{1\leq k\leq k_n,n\in\NN}$
follows from~\eqref{cond:TMa} by dominated convergence. Condition~\eqref{cond:N} 
follows from Propositions~\ref{Proposition:MimpliesCLB}
and~\ref{Proposition:CLBandRimplyN}.\hfill$\Box$
\medskip

\textit{Proof of Proposition \ref{Proposition:EquivalenceforCLB}.} Statement (ii) follows from (i) because 
$\left|X_{nk,j}\right|\leq\left\lVert X_{nk}\right\rVert$ for all $n,k$ and $j$ so that
\begin{displaymath}
\sum_{k=1}^{k_n}E\left(X_{nk,j}^21_{\{\left|X_{nk,j}\right|\geq\varepsilon\}}|\F_{n,k-1}\right)\leq
\sum_{k=1}^{k_n}E\left(\left\lVert X_{nk}\right\rVert^21_{\{\left\lVert X_{nk}\right\rVert\geq\varepsilon\}}|\F_{n,k-1}\right)
\end{displaymath}
for all $n$ and $j$. To obtain statement (i) from (ii) observe that for all $n$ and $k$
\begin{displaymath}
\left\lVert X_{nk}\right\rVert^21_{\{\left\lVert X_{nk}\right\rVert\geq\varepsilon\}}=
\sum_{j=1}^d X_{nk,j}^21_{\{\sum_{\ell=1}^d X_{nk,\ell}^2\geq\varepsilon^2\}}\leq
\sum_{j,\ell=1}^d X_{nk,j}^21_{\{X_{nk,\ell}^2\geq\varepsilon^2/d\}}
\end{displaymath}
where
\begin{align*}
&X_{nk,j}^21_{\{X_{nk,\ell}^2\geq\varepsilon^2/d\}}=X_{nk,j}^2
1_{\{\left|X_{nk,\ell}\right|\geq\varepsilon/d^{1/2}\}}\\
&\quad=X_{nk,j}^21_{\{\left|X_{nk,\ell}\right|\geq\varepsilon/d^{1/2}\}\cap\{\left|X_{nk,\ell}\right|
\leq\left|X_{nk,j}\right|\}}+
X_{nk,j}^21_{\{\left|X_{nk,\ell}\right|\geq\varepsilon/d^{1/2}\}
\cap\{\left|X_{nk,\ell}\right|>\left|X_{nk,j}\right|\}}\\
&\quad\leq 
X_{nk,j}^21_{\{\left|X_{nk,j}\right|\geq\varepsilon/d^{1/2}\}}+X_{nk,\ell}^21_{\{\left|X_{nk,\ell}
\right|\geq\varepsilon/d^{1/2}\}}\,.
\end{align*}
Hence for all $n,k$ and $\varepsilon$
\begin{align*}
&\left\lVert X_{nk}\right\rVert^21_{\{\left\lVert X_{nk}\right\rVert\geq\varepsilon\}}\leq
\sum_{j,\ell=1}^d X_{nk,j}^21_{\{\left|X_{nk,j}\right|\geq\varepsilon/d^{1/2}\}}+
\sum_{j,\ell=1}^d X_{nk,\ell}^21_{\{\left|X_{nk,\ell}\right|\geq\varepsilon/d^{1/2}\}}\\
&\quad=2d\sum_{j=1}^d X_{nk,j}^21_{\{\left|X_{nk,j}\right|\geq\varepsilon/d^{1/2}\}}\,.
\end{align*}
Thus for all $n$ and $\varepsilon$
\begin{displaymath}
\sum_{k=1}^{k_n}E\left(\left\lVert X_{nk}\right\rVert^21_{\{\left\lVert X_{nk}\right\rVert\geq\varepsilon\}}|\F_{n,k-1}\right)\leq
2d\sum_{j=1}^d\sum_{k=1}^{k_n}E\left(X_{nk,j}^2
1_{\{\left|X_{nk,j}\right|\geq\varepsilon/d^{1/2}\}}|\F_{n,k-1}\right)
\end{displaymath}
so that (i) holds. The proofs of the statements about \eqref{cond:CLB1}, \eqref{cond:M1} and \eqref{cond:M2} are
similar.\hfill$\Box$
\medskip

\section{\large Application: Stable autoregressive models of order $\boldsymbol{d}$}
\label{Section:Application}

We consider an autoregressive process $\left(Y_n\right)_{n\geq-d+1}$ of order $d$ generated recursively by
\begin{displaymath}
Y_n=\theta_1Y_{n-1}+\cdots+\theta_dY_{n-d}+Z_n\,,\quad n\geq1\,,
\end{displaymath}
where $\theta=\left(\theta_1,\ldots,\theta_d\right)^T\in\RR^d,\,\left(Z_n\right)_{n\geq1}$ is an i.i.d.
sequence of real random variables with $Z_1\in\LL^2\left(P\right),\,E\left(Z_1\right)=0,\,\sigma^2=\Var\left(Z_1\right)>0$
and the initial value $\left(Y_0,\ldots,Y_{-d+1}\right)^T$ is an $\RR^d$-valued random vector
independent of $\left(Z_n\right)_{n\geq1}$ with $Y_0,\ldots,Y_{-d+1}\in\LL^2\left(P\right)$. The least squares estimator
for the parameter $\theta$ on the basis of the observations $Y_{-d+1},\ldots,Y_0,Y_1,\ldots,Y_n$ is given by
\begin{displaymath}
\widehat{\theta}_n=\left(\sum_{k=1}^nU_{k-1}U_{k-1}^T\right)^{-1}\sum_{k=1}^nU_{k-1}Y_k\,,\quad n\geq1\,,
\end{displaymath}
where $U_k=\left(Y_k,Y_{k-1},\ldots,Y_{k-d+1}\right)^T,k\geq0$, and where 
$\left(\sum_{k=1}^nU_{k-1}U_{k-1}^T\right)^{-1}:=0$ if the random symmetric and positive semidefinite matrix
$\sum_{k=1}^nU_{k-1}U_{k-1}^T$ is not positive definite and therefore singular.
\smallskip

We prove a mixing central limit theorem for $\widehat{\theta}_n$ in the so-called stable case of the autoregressive model. 
For classical convergence in distribution the asymptotic theory for $\widehat{\theta}_n$ is well-known in 
this case; see e.g. \cite{Anderson2}, Chapter~5, \cite{BrockwellDavis}, Chapter~8 and \cite{MannWald}.
\smallskip

The process $\left(U_n\right)_{n\geq0}$ can be expressed as a $d$-dimensional autoregession of order one,
namely
\begin{displaymath}
U_n=BU_{n-1}+W_n\,,\quad n\geq1\,,
\end{displaymath}
where
\begin{displaymath}
B=B\left(\theta\right)=
\begin{pmatrix}
\theta_1 & \ldots  & \theta_{d-1} & \theta_d \\
         &         &              & 0        \\
         & I_{d-1} &              & \vdots   \\
         &         &              & 0        
\end{pmatrix}
\qquad\text{and}\qquad W_n=\left(Z_n,0,\ldots,0\right)^T\,.
\end{displaymath}
Then, by induction, we have
\begin{equation}
\label{eq:Un}
U_n=B^nU_0+\sum_{j=1}^nB^{n-j}W_j=B^nU_0+\sum_{j=1}^nB^{j-1}W_{n+1-j}\,,\quad n\geq0\,.
\end{equation}
Let $\F_n=\sigma\left(U_0,Z_1,\ldots,Z_n\right),n\geq0$. Then all processes are
adapted to the filtration $\left(\F_n\right)_{n\geq0}$.
\smallskip

In the sequel we assume that the spectral radius of $B$ is strictly less than one, i.e. 
all eigenvalues of $B$ are strictly less than one in absolute value.
This condition defines the \textit{stable case} in the theory of autoregression models. Then there
exists a submultiplicative norm $\left\lVert\cdot\right\rVert_B$ on the space $\RR^{d\times d}$ of all real 
$d\times d$-matrices with $\left\lVert B\right\rVert_B<1$; see \cite{Schott}, Theorem 4.24. The Frobenius norm  
$\left\lVert\cdot\right\rVert_F$ on $\RR^{m\times n}$ is also submultiplicative,
satisfies $\left\lVert D\right\rVert_F=\left\lVert D^T\right\rVert_F$ and is compatible with the 
Euclidean $\ell^2$-norm $\left\lVert\cdot\right\rVert$ (same notation on various spaces), i.e. 
$\left\lVert Du\right\rVert\leq\left\lVert D\right\rVert_F\left\lVert u\right\rVert$ 
for all $u\in\RR^n$ and $D\in\RR^{m\times n}$ as a consequence of submultiplicativity and 
$\left\lVert u\right\rVert_F=\left\lVert u\right\rVert$.
Since all norms on $\RR^{d\times d}$ are equivalent, there is a real constant $C\left(B\right)$ such that
$\left\lVert D\right\rVert_F\leq C\left(B\right)\left\lVert D\right\rVert_B$ for all $D\in\RR^{d\times d}$.
\smallskip

Note that the above condition on the eigenvalues of $B$ implies the existence of a stationary distribution $\mu$,
$U_n\stackrel{d}{\rightarrow}\mu$ as $n\to\infty$, where $\mu$ is the distribution of $\sum_{j=0}^\infty B^jW_{j+1}$,
and for the choice $P^{U_0}=\mu$ the process $\left(Y_{n-d+1}\right)_{n\geq0}$ is strictly stationary.
The covariance matrix of $\mu$ is given by
\begin{displaymath}
\sigma^2\sum_{j=0}^\infty B^j\widetilde{I}_d\left(B^j\right)^T\,,
\end{displaymath}
where $\widetilde{I}_d$ denotes the $d\times d$-matrix with a one in the upper left corner and zeros elsewhere.
To see that both series converge observe that $\left\lVert B^j\right\rVert_B\leq\left\lVert B\right\rVert_B^j$
for all $j\geq1$ by submultiplicativity of $\left\lVert\cdot\right\rVert_B$. Consequently,
\begin{align*}
\kappa_1&:=\sum_{j=0}^\infty\left\lVert B^j\right\rVert_F=\left\lVert I_d\right\rVert_F+\sum_{j=1}^\infty\left\lVert B^j\right\rVert_F
\leq d^{1/2}+C\left(B\right)\sum_{j=1}^\infty\left\lVert B^j\right\rVert_B\\
&\leq d^{1/2}+C\left(B\right)\sum_{j=1}^\infty\left\lVert B\right\rVert_B^j<\infty
\end{align*}
because $\left\lVert B\right\rVert_B<1$. Therefore
\begin{align*}
E\left(\sum_{j=0}^\infty\left\lVert B^j\right\rVert_F\left\lvert Z_{j+1}\right\rvert\right)
=\sum_{j=0}^\infty\left\lVert B^j\right\rVert_FE\left(\left\lvert Z_{j+1}\right\rvert\right)
=\kappa_1E\left(\lvert Z_1\rvert\right)<\infty
\end{align*}
so that
\begin{equation}
\label{eq:sumfinite}
\sum_{j=0}^\infty\left\lVert B^jW_{j+1}\right\rVert\leq\sum_{j=0}^\infty\left\lVert B^j\right\rVert_F\left\lVert W_{j+1}\right\rVert
=\sum_{j=0}^\infty\left\lVert B^j\right\rVert_F\left\lvert Z_{j+1}\right\rvert<\infty
\end{equation}
almost surely (which is true already under $E\log^+\left\lvert Z_1\right\rvert<\infty$, cf. \cite{HaeuslerLuschgy},
Lem\-ma~8.1) and
\begin{displaymath}
\sum_{j=0}^\infty\left\lVert B^j\widetilde{I}_d\left(B^j\right)^T\right\rVert_F\leq
\sum_{j=0}^\infty\left\lVert B^j\right\rVert_F\left\lVert\widetilde{I}_d\right\rVert_F\left\lVert\left(B^j\right)^T\right\rVert_F
=\sum_{j=0}^\infty\left\lVert B^j\right\rVert_F^2=:\kappa_2<\infty
\end{displaymath}
in view of $\kappa_1<\infty$. In the following theorem it will be crucial that the symmetric matrix
\begin{displaymath}
\Sigma=\Sigma\left(\vartheta\right)=\sum_{j=0}^\infty B^j\widetilde{I}_d\left(B^j\right)^T
\end{displaymath}
is positive definite and therefore invertible. This will be shown at the end of this section in Lemma~\ref{Lemma:Invertibility}.
The $\sigma$-field $\F_\infty=\sigma\left(\bigcup_{n=0}^\infty\F_n\right)$ 
is defined by the filtration $\left(\F_n\right)_{n\geq0}$. 
\medskip

\begin{Theorem}
\label{Theorem:CLTsforAR}
In the above setting,
\begin{displaymath}
\sqrt{n}\left(\widehat{\theta}_n-\theta\right)\rightarrow\Sigma^{-1/2}N_d\quad\F_\infty\text{-mixing}
\end{displaymath}
and
\begin{displaymath}
\left(\sum_{k=1}^nU_{k-1}U_{k-1}^T\right)^{1/2}
\left(\widehat{\theta}_n-\theta\right)\rightarrow\sigma N_d\quad\F_\infty\text{-mixing}
\end{displaymath}
as $n\to\infty$, where $P^{N_d}=N\left(0,I_d\right)$, $N_d$ is independent of $\F_\infty$ and
$\Sigma^{-1/2}:=\left(\Sigma^{1/2}\right)^{-1}$.
\end{Theorem}
\medskip

The above statements may also be read as 
\begin{displaymath}
\sqrt{n}\left(\widehat{\theta}_n-\theta\right)\rightarrow N\left(0,\Sigma^{-1}\right)\quad\text{mixing}
\end{displaymath}
and
\begin{displaymath}
\left(\sum_{k=1}^nU_{k-1}U_{k-1}^T\right)^{1/2}
\left(\widehat{\theta}_n-\theta\right)\rightarrow N\left(0,\sigma^2I_d\right)\quad\text{mixing}\,,
\end{displaymath}
where mixing is short for $\F$-mixing. In case $d=1$ Theorem~\ref{Theorem:CLTsforAR} reduces to
The\-o\-rem~9.1 in \cite{HaeuslerLuschgy}.
\smallskip

We will present here two proofs of Theorem~\ref{Theorem:CLTsforAR}. Both proofs are based on the fact
that the process $\left(X_n\right)_{n\geq1}$ with $X_n=U_{n-1}Z_n$ is a square integrable martingale
difference sequence. In the first proof we use the filtration $\left(\F_n\right)_{n\geq0}$ introduced above and apply the
stable central limit theorem of Corollary~\ref{Corollary:Lindebergforsequences} with $K_n=n^{-1/2}I_d$,
and stationarity plays no role. 
The conditional Lindeberg condition~\eqref{eq:CLBforsequences} and the norming 
condition~\eqref{eq:normingforsequences} will be verified by elementary but lengthy computations exploiting nothing
more than the dynamics of the $U_n$-process via the representation~\eqref{eq:Un}. The second proof uses the
fact that for a starting value $\U_0$ with $P^{\U_0}=\mu$ the process $\left(X_n\right)_{n\geq1}$
is strictly stationary and ergodic and is a martingale difference sequence w.r.t. a filtration
$\left(\G_n\right)_{n\geq0}$ defined below so
that Corollary~\ref{Corollary:Stationarymartingaledifferences} can be applied in this special case. In a
second step it is not difficult to show that the asymptotics of the normalized sums $n^{-1/2}\sum_{k=1}^nX_k$
is the same for all starting values.
\medskip 

\textit{First proof of Theorem~\ref{Theorem:CLTsforAR}.} As announced above, we apply the stable 
central limit theorem of Corollary \ref{Corollary:Lindebergforsequences} with $K_n=n^{-1/2}I_d$
to the square integrable martingale difference sequence $\left(X_n\right)_{n\geq1}$ with $X_n=U_{n-1}Z_n$ 
w.r.t. the filtration $\left(\F_n\right)_{n\geq0}$.
\smallskip

\textit{Step 1.} In a first step, we prove the conditional Lindeberg condition 
\eqref{eq:CLBforsequences}. For all $n\geq1$ and $\varepsilon,\delta>0$ we have
\begin{align*}
&\sum_{k=1}^nE\left(\left\lVert K_nX_k\right\rVert^21_{\{\left\lVert K_nX_k\right\rVert\geq\varepsilon\}}|\F_{k-1}\right)
=\frac{1}{n}\sum_{k=1}^nE\left(\left\lVert U_{k-1}Z_k\right\rVert^21_{\{\left\lVert U_{k-1}Z_k\right\rVert\geq\varepsilon\sqrt{n}\}}|\F_{k-1}\right)\\
&\quad=\frac{1}{n}\sum_{k=1}^n\left\lVert U_{k-1}\right\rVert^2E\left(Z_k^21_{\{\left\lVert U_{k-1}Z_k\right\rVert\geq\varepsilon\sqrt{n}\}
\cap\{\left\lVert U_{k-1}\right\rVert\leq\delta\sqrt{n}\}}|\F_{k-1}\right)\\
&\quad\quad+\frac{1}{n}\sum_{k=1}^n\left\lVert U_{k-1}\right\rVert^2E\left(Z_k^21_{\{\left\lVert U_{k-1}Z_k\right\rVert\geq\varepsilon\sqrt{n}\}
\cap\{\left\lVert U_{k-1}\right\rVert>\delta\sqrt{n}\}}|\F_{k-1}\right)\\
&\quad\leq\left(\frac{1}{n}\sum_{k=1}^n\left\lVert U_{k-1}\right\rVert^2\right)E\left(Z_1^21_{\{\left\lvert Z_1\right\rvert\geq\varepsilon/\delta\}}\right)+
\frac{1}{n}\sum_{k=1}^n\left\lVert U_{k-1}\right\rVert^21_{\{\left\lVert U_{k-1}\right\rVert>\delta\sqrt{n}\}}\sigma^2\\
&\quad=:I_n\left(\varepsilon,\delta\right)+\II_n\left(\delta\right)\sigma^2\,.
\end{align*}
For all $\delta>0$ we will show that $\II_n\left(\delta\right)\rightarrow0$ in probability as $n\to\infty$. For this,
note that for all $n\geq1$ and $\delta,\widetilde{\varepsilon}>0$ we have
\begin{displaymath}
P\left(\frac{1}{n}\sum_{k=1}^n\left\lVert U_{k-1}\right\rVert^21_{\{\left\lVert U_{k-1}\right\rVert>\delta\sqrt{n}\}}\geq\widetilde{\varepsilon}\right)
\leq P\left(\max_{1\leq k\leq n}\left\lVert U_{k-1}\right\rVert>\delta\sqrt{n}\right)
\end{displaymath}
so that it is enough to show that
\begin{displaymath}
\frac{1}{\sqrt{n}}\max_{1\leq k\leq n}\left\lVert U_{k-1}\right\rVert\rightarrow0\quad\text{in probability as }n\to\infty\,.
\end{displaymath}
From \eqref{eq:Un} we obtain for all $k\geq1$
\begin{align*}
\left\lVert U_{k-1}\right\rVert&\leq\left\lVert B^{k-1}U_0\right\rVert+\sum_{j=0}^{k-2}\left\lVert B^jW_{k-1-j}\right\rVert\leq
\left\lVert B^{k-1}\right\rVert_F\left\lVert U_0\right\rVert+\sum_{j=0}^{k-2}\left\lVert B^j\right\rVert_F\left\lVert W_{k-1-j}\right\rVert\\
&\leq\kappa_3\left\lVert U_0\right\rVert+\sum_{j=0}^{k-2}\left\lVert B^j\right\rVert_F\left\lvert Z_{k-1-j}\right\rvert
\leq\kappa_3\left\lVert U_0\right\rVert+\kappa_1\max_{1\leq j\leq k-1}\left\lvert Z_j\right\rvert\,,
\end{align*}
where $\kappa_3:=\sup_{j\geq 0}\left\lVert B^j\right\rVert_F<\infty$ because of $\kappa_1<\infty$.
Consequently, for all $n\geq1$,
\begin{displaymath}
\frac{1}{\sqrt{n}}\max_{1\leq k\leq n}\left\lVert U_{k-1}\right\rVert\leq\frac{1}{\sqrt{n}}
\kappa_3\left\lVert U_0\right\rVert+\kappa_1\frac{1}{\sqrt{n}}\max_{1\leq j\leq n}\left\lvert Z_j\right\rvert\,,
\end{displaymath}
where the right-hand side converges in probability to zero as $n\to\infty$ because $\left(Z_i\right)_{i\geq1}$ is an
i.i.d. sequence with $E\left(Z_1^2\right)<\infty$, and $\II_n\left(\delta\right)\rightarrow0$ in probability as $n\to\infty$  is proven.
From the above bound on $\left\lVert U_{k-1}\right\rVert$ we also obtain for all $k\geq1$
\begin{align*}
E\left(\left\lVert U_{k-1}\right\rVert^2\right)^{1/2}& \leq\kappa_3E\left(\left\lVert U_0\right\rVert^2\right)^{1/2}
+\sum_{j=0}^{k-2}\left\lVert B^j\right\rVert_F E\left(Z_1^2\right)^{1/2}\\
&\leq\kappa_3E\left(\left\lVert U_0\right\rVert^2\right)^{1/2}+\kappa_1\sigma=:\kappa_4\,.
\end{align*}
This gives $E\left(\left\lVert U_{k-1}\right\rVert^2\right)\leq\kappa_4^2$ for all $k\geq1$ so that for all $n\geq1$ and
$\varepsilon,\delta>0$
\begin{displaymath}
E\left(\left\lvert I_n\left(\varepsilon,\delta\right)\right\rvert \right)=E\left(\frac{1}{n}\sum_{k=1}^n\left\lVert U_{k-1}\right\rVert^2\right)
E\left(Z_1^21_{\{\left\lvert Z_1\right\rvert>\varepsilon/\delta\}}\right)\leq\kappa_4^2
E\left(Z_1^21_{\{\left\lvert Z_1\right\rvert>\varepsilon/\delta\}}\right)\,.
\end{displaymath}
Hence for all $n\geq1$ and $\varepsilon,\widetilde{\varepsilon},\delta>0$
\begin{align*}
&P\left(\left\lvert I_n\left(\varepsilon,\delta\right)+\II_n\left(\delta\right)\sigma^2\right\rvert\geq\widetilde{\varepsilon}\right)
\leq P\left(\left\lvert I_n\left(\varepsilon,\delta\right)\right\rvert\geq\widetilde{\varepsilon}/2\right)+
P\left(\left\lvert\II_n\left(\delta\right)\right\rvert\geq\widetilde{\varepsilon}/2\sigma^2\right)\\
&\quad\leq\frac{2}{\widetilde{\varepsilon}}E\left(\left\lvert I_n\left(\varepsilon,\delta\right)\right\rvert\right)+
P\left(\left\lvert\II_n\left(\delta\right)\right\rvert\geq\widetilde{\varepsilon}/2\sigma^2\right)\\
&\quad\leq\frac{2\kappa_4^2}{\widetilde{\varepsilon}}E\left(Z_1^21_{\{\left\lvert Z_1\right\rvert>\varepsilon/\delta\}}\right)+
P\left(\left\lvert\II_n\left(\delta\right)\right\rvert\geq\widetilde{\varepsilon}/2\sigma^2\right)\,.
\end{align*}
Note that $\varepsilon,\widetilde{\varepsilon}>0$ are fixed. Because of $\II_n\left(\delta\right)\rightarrow0$ in probability
as $n\to\infty$ and $E\left(Z_1^2\right)<\infty$ the right-hand
side of this inequality can be made arbitrarily small for all large $n$ by choosing $\delta$ sufficiently small.
This concludes the proof of \eqref{eq:CLBforsequences}.
\smallskip

\textit{Step 2.} We now turn to the proof of \eqref{eq:normingforsequences} and show
\begin{equation}
\label{eq:NinAR}
K_n\sum_{k=1}^nE\left(X_kX_k^T|\F_{k-1}\right)K_n^T\rightarrow\sigma^4\Sigma\quad\text{in probability as }n\to\infty\,.
\end{equation}
For all $n\geq1$ we have
\begin{displaymath}
K_n\sum_{k=1}^nE\left(X_kX_k^T|\F_{k-1}\right)K_n^T=\frac{1}{n}
\sum_{k=1}^nE\left(Z_kU_{k-1}U_{k-1}^TZ_k|\F_{k-1}\right)=\frac{\sigma^2}{n}\sum_{k=1}^nU_{k-1}U_{k-1}^T
\end{displaymath}
so that \eqref{eq:NinAR} follows from
\begin{equation}
\label{eq:convergenceofUUT}
\frac{1}{n}\sum_{k=1}^nU_{k-1}U_{k-1}^T\rightarrow\sigma^2\Sigma\quad\text{in probability as }n\to\infty\,.
\end{equation}
For all $n\geq1$ we obtain from \eqref{eq:Un}
\begin{align*}
&\frac{1}{n}\sum_{k=1}^nU_{k-1}U_{k-1}^T\\
&\quad=\frac{1}{n}\sum_{k=1}^nB^{k-1}U_0U_0^T\left(B^{k-1}\right)^T+
\frac{1}{n}\sum_{k=1}^nB^{k-1}U_0\left(\sum_{j=1}^{k-1}B^{j-1}W_{k-j}\right)^T\\
&\qquad+\frac{1}{n}\sum_{k=1}^n\left(\sum_{j=1}^{k-1}B^{j-1}W_{k-j}\right)U_0^T\left(B^{k-1}\right)^T\\
&\qquad+\frac{1}{n}\sum_{k=1}^n\left(\sum_{j=1}^{k-1}B^{j-1}W_{k-j}\right)\left(\sum_{j=1}^{k-1}B^{j-1}W_{k-j}\right)^T\\
&\quad=:I_{1,n}+I_{2,n}+I_{3,n}+I_{4,n}\,.
\end{align*}
For all $n\geq1$, let $G_{1,n}\left(I\right)=B^{n-1}U_0U_0^T\left(B^{n-1}\right)^T$. Then
\begin{align*}
\left\lVert G_{1,n}\left(I\right)\right\rVert_F\leq\left\lVert B^{n-1}\right\rVert_F
\left\lVert U_0\right\rVert^2\left\lVert B^{n-1}\right\rVert_F
\leq\left\lVert B^{n-1}\right\rVert_F^2\left\lVert U_0\right\rVert^2
\rightarrow0
\end{align*}
everywhere because $\kappa_1<\infty$ so that
\begin{displaymath}
\left\lVert I_{1,n}\right\rVert_F=\left\lVert\frac{1}{n}\sum_{k=1}^nG_{1,k}\left(I\right)\right\rVert_F
\leq\frac{1}{n}\sum_{k=1}^n\left\lVert G_{1,k}\left(I\right)\right\rVert_F\rightarrow0\quad\text{everywhere}
\end{displaymath}
which implies $I_{1,n}\rightarrow0$ in probability as $n\to\infty$. 
\smallskip

Our next aim is to show that
$I_{2,n}\rightarrow0$ in probability as $n\to\infty$. For $n\geq1$, let
$G_{2,n}\left(I\right)=B^{n-1}U_0\left(\sum_{j=1}^{n-1}B^{j-1}W_{n-j}\right)^T$. Then
\begin{align*}
\left\lVert G_{2,n}\left(I\right)\right\rVert_F&\leq\left\lVert B^{n-1}\right\rVert_F\left\lVert U_0\right\rVert
\left\lVert\sum_{j=1}^{n-1}B^{j-1}W_{n-j}\right\rVert\\
&\leq\left\lVert B^{n-1}\right\rVert_F\left\lVert U_0\right\rVert
\sum_{j=1}^{n-1}\left\lVert B^{j-1}\right\rVert_F\left\lvert Z_{n-j}\right\rvert
\end{align*}
so that
\begin{align*}
E\left(\left\lVert G_{2,n}\left(I\right)\right\rVert_F\right)&\leq\left\lVert B^{n-1}\right\rVert_F
E\left(\left\lVert U_0\right\rVert\right)E\left(\left\lvert Z_1\right\rvert\right)
\sum_{j=1}^{n-1}\left\lVert B^{j-1}\right\rVert_F\\
&\leq\left\lVert B^{n-1}\right\rVert_FE\left(\left\lVert U_0\right\rVert\right)E\left(\left\lvert Z_1\right\rvert\right)
\kappa_1\rightarrow0\,
\end{align*}
as $n\to\infty$ because $\kappa_1<\infty$. Consequently,
\begin{displaymath}
E\left(\left\lVert I_{2,n}\right\rVert_F\right)=E\left(\left\lVert\frac{1}{n}\sum_{k=1}^nG_{2,k}\left(I\right)\right\rVert_F\right)
\leq\frac{1}{n}\sum_{k=1}^nE\left(\left\lVert G_{2,k}\left(I\right)\right\rVert_F\right)\rightarrow0\quad\text{as }n\to\infty
\end{displaymath}
which proves $I_{2,n}\rightarrow0$ in probability as $n\to\infty$. Because of $I_{3,n}=I_{2,n}^T$ we also have
$I_{3,n}\rightarrow0$ in probability as $n\to\infty$. Therefore, it remains to consider $I_{4,n}$. For all $n\geq1$
we have
\begin{align*}
I_{4,n}&=\frac{1}{n}\sum_{k=1}^n\sum_{j,m=1}^{k-1}B^{j-1}W_{k-j}W_{k-m}^T\left(B^{m-1}\right)^T\\
&=\frac{1}{n}\sum_{k=1}^n\sum_{j,m=1}^{k-1}B^{k-j-1}W_{j}W_{m}^T\left(B^{k-m-1}\right)^T\\
&=\frac{1}{n}\sum_{k=1}^n\sum_{j,m=1}^{k-1}B^{k-j-1}\widetilde{I}_d\left(B^{k-m-1}\right)^TZ_jZ_m\\
&=\frac{1}{n}\sum_{k=1}^n\sum_{j=1}^{k-1}B^{k-j-1}\widetilde{I}_d\left(B^{k-j-1}\right)^TZ_j^2\\
&\quad+\frac{1}{n}\sum_{k=1}^n\sum_{1\leq j<m\leq k-1}B^{k-j-1}\widetilde{I}_d\left(B^{k-m-1}\right)^TZ_jZ_m\\
&\quad+\frac{1}{n}\sum_{k=1}^n\sum_{1\leq m<j\leq k-1}B^{k-j-1}\widetilde{I}_d\left(B^{k-m-1}\right)^TZ_jZ_m
=:\II_{1,n}+\II_{2,n}+\II_{3,n}\,.
\end{align*}
We will examine $\II_{3,n}$ first. For all $n\geq1$ we have
\begin{align*}
&\II_{3,n}=\frac{1}{n}\sum_{k=1}^n\sum_{j=2}^{k-1}\sum_{m=1}^{j-1}B^{k-1-j}\widetilde{I}_d\left(B^{k-1-m}\right)^TZ_jZ_m\\
&= \frac{1}{n}\sum_{j=2}^{n-1}\sum_{k=j+1}^n\sum_{m=1}^{j-1}B^{k-1-j}\widetilde{I}_d\left(B^{k-1-m}\right)^TZ_jZ_m\\
&=\frac{1}{n}\sum_{j=2}^{n-1}\sum_{m=1}^{j-1}\left(\sum_{k=j+1}^nB^{k-j-1}\widetilde{I}_d\left(B^{k-m-1}\right)^T\right)Z_jZ_m
=:\frac{1}{n}\sum_{j=2}^{n-1}\sum_{m=1}^{j-1}M_n\left(j,m\right)Z_jZ_m\,,
\end{align*}
with deterministic $d\times d$-matrices $M_n\left(j,m\right)$. Now for all $n\geq2$
\begin{displaymath}
E\left(\left\lVert\II_{3,n}\right\rVert_F^2\right)=\frac{1}{n^2}\sum_{p,q=1}^dE\left(\left(\sum_{j=2}^{n-1}\left[\sum_{m=1}^{j-1}
M_n\left(j,m\right)_{p,q}Z_m\right]Z_j\right)^2\right)\,.
\end{displaymath}
Since the $Z_j$ are independent with mean zero we obtain for $1\leq j^\prime<j$ that
\begin{align*}
&E\left(\left[\sum_{m=1}^{j^\prime-1}M_n\left(j^\prime,m\right)_{p,q}Z_m\right]Z_{j^\prime}
\left[\sum_{m=1}^{j-1}M_n\left(j,m\right)_{p,q}Z_m\right]Z_j\right)\\
&\quad=E\left(\left[\sum_{m=1}^{j^\prime-1}M_n\left(j^\prime,m\right)_{p,q}Z_m\right]Z_{j^\prime}
\left[\sum_{m=1}^{j-1}M_n\left(j,m\right)_{p,q}Z_m\right]\right)E\left(Z_j\right)=0
\end{align*}
and likewise for $1\leq j<j^\prime$. Therefore we get
\begin{align*}
&E\left(\left\lVert\II_{3,n}\right\rVert_F^2\right)=\frac{1}{n^2}\sum_{p,q=1}^d\sum_{j=2}^{n-1}E\left(
\left[\sum_{m=1}^{j-1}M_n\left(j,m\right)_{p,q}Z_m\right]^2Z_j^2\right)\\
&\quad=\frac{\sigma^2}{n^2}\sum_{p,q=1}^d\sum_{j=2}^{n-1}E\left(
\left[\sum_{m=1}^{j-1}M_n\left(j,m\right)_{p,q}Z_m\right]^2\right)\\
&\quad=\frac{\sigma^2}{n^2}\sum_{p,q=1}^d\sum_{j=2}^{n-1}\sum_{m=1}^{j-1}E\left(\left[M_n\left(j,m\right)_{p,q}Z_m\right]^2\right)\\
&\quad=\frac{\sigma^2}{n^2}\sum_{p,q=1}^d\sum_{j=2}^{n-1}\sum_{m=1}^{j-1}M_n\left(j,m\right)_{p,q}^2E\left(Z_m^2\right)
=\frac{\sigma^4}{n^2}\sum_{j=2}^{n-1}\sum_{m=1}^{j-1}\left\lVert M_n\left(j,m\right)\right\rVert_F^2\,,
\end{align*}
where the next to the last equation follows as before from the independence of the $Z_j$ and $E\left(Z_j\right)=0$.
Now for all $2\leq j\leq n-1$ and $1\leq m\leq j-1$
\begin{align*}
&\left\lVert M_n\left(j,m\right)\right\rVert_F\leq\sum_{k=j+1}^n\left\lVert B^{k-j-1}
\right\rVert_F\left\lVert\widetilde{I}_d\right\rVert_F\left\lVert\left(B^{k-m-1}\right)^T\right\rVert_F\\
&\quad=\sum_{k=j+1}^n\left\lVert B^{k-j-1}\right\rVert_F\left\lVert B^{k-m-1}\right\rVert_F
=\sum_{\ell=0}^{n-j-1}\left\lVert B^\ell\right\rVert_F\left\lVert B^{\ell+j-m}\right\rVert_F\\ 
&\quad\leq\left\lVert B^{j-m}\right\rVert_F\sum_{\ell=0}^\infty\left\lVert B^\ell\right\rVert_F^2
=\kappa_2\left\lVert B^{j-m}\right\rVert_F\,.
\end{align*}
Using this bound we obtain for all $n\geq1$
\begin{displaymath}
E\left(\left\lVert\II_{3,n}\right\rVert_F^2\right)\leq\sigma^4\kappa_2^2\frac{1}{n^2}
\sum_{j=2}^{n-1}\sum_{m=1}^{j-1}\left\lVert B^{j-m}\right\rVert_F^2\leq
\sigma^4\kappa_2^3\frac{1}{n}
\end{displaymath}
which proves $\II_{3,n}\rightarrow0$ in probability as $n\to\infty$. Because of $\II_{2,n}=\II_{3,n}^T$ we also have
$\II_{2,n}\rightarrow0$ in probability as $n\to\infty$. It remains to consider
\begin{displaymath}
\II_{1,n}=\frac{1}{n}\sum_{k=1}^n\sum_{j=0}^{k-2}B^j\widetilde{I}_d\left(B^j\right)^TZ_{k-j-1}^2\,.
\end{displaymath}
For this, we write
\begin{align*}
\II_{1,n}-\sigma^2\Sigma=&\II_{1,n}-\sigma^2\sum_{j=0}^\infty B^j\widetilde{I}_d\left(B^j\right)^T
=\frac{1}{n}\sum_{k=1}^n\sum_{j=0}^{k-2}B^j\widetilde{I}_d\left(B^j\right)^T\left(Z_{k-j-1}^2-\sigma^2\right)\\
&\quad-\frac{\sigma^2}{n}\sum_{k=1}^n\sum_{j=k-1}^\infty B^j\widetilde{I}_d\left(B^j\right)^T
=:\III_{1,n}-\sigma^2\III_{2,n}\,.
\end{align*}
As for the deterministic $\III_{2,n}$, let, for $n\geq1$, $G_{2,n}\left(\III\right)=
\sum_{j=n-1}^\infty B^j\widetilde{I}_d\left(B^j\right)^T$. Then
\begin{displaymath}
\left\lVert G_{2,n}\left(\III\right)\right\rVert_F\leq\sum_{j=n-1}^\infty\left\lVert B^j\widetilde{I}_d
\left(B^j\right)^T\right\rVert_F\rightarrow0\quad\text{as }n\to\infty
\end{displaymath}
because $\kappa_2<\infty$ so that
\begin{displaymath}
\left\lVert\III_{2,n}\right\rVert_F=\left\lVert\frac{1}{n}\sum_{k=1}^nG_{2,k}\left(\III\right)\right\rVert_F
\leq\frac{1}{n}\sum_{k=1}^n\left\lVert G_{2,k}\left(\III\right)\right\rVert_F\rightarrow0\quad\text{as }n\to\infty\,.
\end{displaymath}
It remains to show that $\III_{1,n}\rightarrow0$ in probability as $n\to\infty$. For this, we set
$V_j=Z_j^2-\sigma^2$ for all $j\geq1$. These random variables are i.i.d. and integrable with mean zero. For all $n\geq1$ we have
\begin{align*}
&\III_{1,n}=\frac{1}{n}\sum_{k=1}^n\sum_{j=0}^{k-2}B^j\widetilde{I}_d\left(B^j\right)^TV_{k-j-1}=\frac{1}{n}
\sum_{k=1}^n\sum_{j=1}^{k-1}B^{k-j-1}\widetilde{I}_d\left(B^{k-j-1}\right)^TV_j\\
&=\frac{1}{n}\sum_{j=1}^{n-1}\sum_{k=j+1}^nB^{k-j-1}\widetilde{I}_d\left(B^{k-j-1}\right)^TV_j
=:\frac{1}{n}\sum_{j=1}^{n-1}R_n\left(j\right)V_j\,,
\end{align*}
with deterministic $d\times d$-matrices $R_n\left(j\right)$. Because of $E\left(V_j\right)=0$ for all $j\geq1$,
for all $n\geq1$ and $0<c<\infty$ we can write
\begin{align*}
\III_{1,n}&=\frac{1}{n}\sum_{j=1}^{n-1}R_n\left(j\right)\left[V_j1_{\{\left\lvert V_j\right\rvert\leq c\}}
-E\left(V_j1_{\{\left\lvert V_j\right\rvert\leq c\}}\right)\right]\\
&\ \ +\frac{1}{n}\sum_{j=1}^{n-1}R_n\left(j\right)V_j1_{\{\left\lvert V_j\right\rvert>c\}}-
\frac{1}{n}\sum_{j=1}^{n-1}R_n\left(j\right)E\left(V_j1_{\{\left\lvert V_j\right\rvert>c\}}\right)\\
&=:\IV_{1,n}+\IV_{2,n}-\IV_{3,n}\,.
\end{align*}
For all $n\geq2$ and $1\leq j\leq n-1$ we have
\begin{align*}
\left\lVert R_n\left(j\right)\right\rVert_F\leq\sum_{k=j+1}^n\left\lVert B^{k-j-1}
\right\rVert_F\left\lVert\widetilde{I}_d\right\rVert_F\left\lVert\left(B^{k-j-1}\right)^T\right\rVert_F
=\sum_{k=j+1}^n\left\lVert B^{k-j-1}\right\rVert_F^2\leq\kappa_2\,.
\end{align*}
This bound yields for all $n\geq1$
\begin{displaymath}
E\left(\left\lVert\IV_{2,n}\right\rVert_F\right)\leq\frac{1}{n}\sum_{j=1}^{n-1}\left\lVert R_n\left(j\right)\right\rVert_F
E\left(\left\lvert V_1\right\rvert1_{\{\left\lvert V_1\right\rvert>c\}}\right)\leq\kappa_2
E\left(\left\lvert V_1\right\rvert1_{\{\left\lvert V_1\right\rvert>c\}}\right)\end{displaymath}
and
\begin{displaymath}
\left\lVert\IV_{3,n}\right\rVert_F\leq\frac{1}{n}\sum_{j=1}^{n-1}\left\lVert R_n\left(j\right)\right\rVert_F
E\left(\left\lvert V_1\right\rvert1_{\{\left\lvert V_1\right\rvert>c\}}\right)\leq\kappa_2
E\left(\left\lvert V_1\right\rvert1_{\{\left\lvert V_1\right\rvert>c\}}\right)\,.
\end{displaymath}
To produce a bound for $\IV_{1,n}$ we set $V_j\left(c\right)=V_j1_{\{\left\lvert V_j\right\rvert\leq c\}}
-E\left(V_j1_{\{\left\lvert V_j\right\rvert\leq c\}}\right)$
and get
\begin{align*}
&E\left(\left\lVert\IV_{1,n}\right\rVert_F^2\right)=\frac{1}{n^2}
E\left(\left\lVert\sum_{j=1}^{n-1}R_n\left(j\right)V_j\left(c\right)\right\rVert_F^2\right)\\
&\quad=\frac{1}{n^2}\sum_{p,q=1}^dE\left(\left(\sum_{j=1}^{n-1}R_n\left(j\right)_{p,q}V_j\left(c\right)\right)^2\right)=
\frac{1}{n^2}\sum_{p,q=1}^d\sum_{j=1}^{n-1}E\left(R_n\left(j\right)_{p,q}^2V_j\left(c\right)^2\right)\,,
\end{align*}
using the fact that the $R_n\left(j\right)_{p,q}$ are deterministic and the random variables $V_j\left(c\right)$ are independent
with mean zero. Because of $\left\lvert V_j\left(c\right)\right\rvert\leq2c$ we obtain
\begin{displaymath}
E\left(\left\lVert\IV_{1,n}\right\rVert_F^2\right)\leq4c^2\frac{1}{n^2}\sum_{j=1}^{n-1}\left\lVert R_n\left(j\right)\right\rVert_F^2
\leq 4c^2\kappa_2^2\frac{1}{n}\,.
\end{displaymath}
Combining the bounds established so far we arrive at
\begin{align*}
&E\left(\left\lVert\III_{1,n}\right\rVert_F\right)\leq E\left(\left\lVert\IV_{1,n}\right\rVert_F\right)+
E\left(\left\lVert\IV_{2,n}\right\rVert_F\right)+\left\lVert\IV_{3,n}\right\rVert_F\\
&\quad\leq E\left(\left\lVert\IV_{1,n}\right\rVert_F^2\right)^{1/2}+E\left(\left\lVert\IV_{2,n}\right\rVert_F\right)+\left\lVert\IV_{3,n}\right\rVert_F\\
&\quad\leq2c\kappa_2\frac{1}{\sqrt{n}}+2\kappa_2E\left(\left\lvert V_1\right\rvert1_{\{\left\lvert V_1\right\rvert>c\}}\right)
\end{align*}
for every $0<c<\infty$.
Because of $E\left(\left\lvert V_1\right\rvert\right)<\infty$ the right-hand side of this inequality can be made arbitrarily
small for all sufficiently large $n$ by choosing $c$ large enough. This proves $E\left(\left\lVert \III_{1,n}\right\rVert_F\right)\rightarrow0$
as $n\to\infty$, and $\III_{1,n}\rightarrow0$ in probability follows. This concludes the proof of
\eqref{eq:convergenceofUUT} and therefore of \eqref{eq:NinAR}.
\smallskip

\textit{Step 3.} Corollary \ref{Corollary:Lindebergforsequences} now implies
\begin{equation}
\label{eq:CLTforARwrongnorming1}
\frac{1}{\sqrt{n}}\sum_{k=1}^nU_{k-1}Z_k=K_n\sum_{k=1}^nX_k\rightarrow\sigma^2\Sigma^{1/2}N_d\quad\F_\infty\text{-mixing as }n\to\infty\,.
\end{equation}
To complete the proof of the theorem we will rely on the following fact: Let $\widetilde{X},\widetilde{X}_n,n\geq1,$ be
$d$-dimensional random vectors defined on some probability space $\left(\Omega,\F,P\right)$, let $\Omega_n\in\F,n\geq1,$ be events
with $P\left(\Omega_n\right)\rightarrow1$ as $n\to\infty$, and let $\G\subset\F$ be a sub-$\sigma$-field of $\F$. Then, as $n\to\infty$,
\begin{equation}
\label{eq:equivalence}
\widetilde{X}_n\rightarrow \widetilde{X}\quad\G\text{-mixing}\qquad\text{if and only if}\qquad
\widetilde{X}_n1_{\Omega_n}\rightarrow \widetilde{X}\quad\G\text{-mixing}\,.
\end{equation}
This follows from Theorem 3.18 (a) in \cite{HaeuslerLuschgy} because for all $n\geq1$ and $\varepsilon>0$ we have
$P\left(\left\lVert\widetilde{X}_n-\widetilde{X}_n1_{\Omega_n}\right\rVert\geq\varepsilon\right)\leq P\left(\Omega_n^c\right)\rightarrow0$
as $n\to\infty$ so that $\left\lVert\widetilde{X}_n-\widetilde{X}_n1_{\Omega_n}\right\rVert\rightarrow0$ in probability.
We will apply \eqref{eq:equivalence} with
$\Omega_n=\left\{\det\left(\sum_{k-1}^nU_{k-1}U_{k-1}^T\right)>0\right\}$ for $n\geq1,$ for which
\begin{displaymath}
P\left(\Omega_n^c\right)\leq P\left(\left\lvert\det\left(\frac{1}{n}\sum_{k=1}^nU_{k-1}U_{k-1}^T\right)-
\det\left(\sigma^2\Sigma\right)\right\rvert>\frac{1}{2}\det\left(\sigma^2\Sigma\right)\right)\rightarrow0
\end{displaymath}
as $n\to\infty$ by \eqref{eq:convergenceofUUT} because the determinant of a matrix is a continuous function of the matrix components. 
Define $f:\RR^{d\times d}\times\RR^d\rightarrow\RR^d$ by
\begin{displaymath}
f\left(D,x\right)=\left\{
\begin{array}{lcl}
D^{-1}x & , & D\in{\rm GL}\left(d,\RR\right)\,,\\
0       & , & \text{otherwise\,.}
\end{array}\right.
\end{displaymath}
Then $f$ is Borel measurable and continuous at every point in the open subset ${\rm GL}\left(d,\RR\right)\times\RR^d$.
We have on $\Omega_n$
\begin{align*}
\widehat{\theta}_n-\theta&=\left(\sum_{k=1}^nU_{k-1}U_{k-1}^T\right)^{-1}\sum_{k=1}^nU_{k-1}Y_k-\theta\\
&=\left(\sum_{k=1}^nU_{k-1}U_{k-1}^T\right)^{-1}\sum_{k=1}^nU_{k-1}\left(U_{k-1}^T\theta+Z_k\right)-\theta\notag\\
&=\left(\sum_{k=1}^nU_{k-1}U_{k-1}^T\right)^{-1}\sum_{k=1}^nU_{k-1}Z_k\notag
\end{align*}
because $\sum_{k=1}^nU_{k-1}U_{k-1}^T$ is invertible so that still on $\Omega_n$
\begin{displaymath}
\sqrt{n}\left(\widehat{\theta}_n-\theta\right)=f\left(\frac{1}{n}\sum_{k=1}^nU_{k-1}U_{k-1}^T,
\frac{1}{\sqrt{n}}\sum_{k=1}^nU_{k-1}Z_k\right)=:C_n
\end{displaymath}
and
\begin{align*}
\left(\sum_{k=1}^nU_{k-1}U_{k-1}^T\right)^{1/2}\left(\widehat{\theta}_n-\theta\right)&=
f\left(\left(\frac{1}{n}\sum_{k=1}^nU_{k-1}U_{k-1}^T\right)^{1/2},\frac{1}{\sqrt{n}}\sum_{k=1}^nU_{k-1}Z_k\right)\\
&=:D_n\,.
\end{align*}
It follows from \eqref{eq:convergenceofUUT}, \eqref{eq:CLTforARwrongnorming1}
and Theorem~3.18 (b) in \cite{HaeuslerLuschgy} that
\begin{displaymath}
\left(\frac{1}{n}\sum_{k=1}^nU_{k-1}U_{k-1}^T,\frac{1}{\sqrt{n}}\sum_{k=1}^nU_{k-1}Z_k\right)\rightarrow
\left(\sigma^2\Sigma,\sigma^2\Sigma^{1/2}N_d\right)\quad\F_\infty\text{-mixing}
\end{displaymath}
and, using the continuity of the square root of positive semidefinite matrices,
\begin{displaymath}
\left(\left(\frac{1}{n}\sum_{k=1}^nU_{k-1}U_{k-1}^T\right)^{1/2},\frac{1}{\sqrt{n}}\sum_{k=1}^nU_{k-1}Z_k\right)\rightarrow
\left(\sigma\Sigma^{1/2},\sigma^2\Sigma^{1/2}N_d\right)\quad\F_\infty\text{-mixing}
\end{displaymath}
as $n\to\infty$. Using Theorem~3.18 (c) in \cite{HaeuslerLuschgy}, this yields
\begin{displaymath}
C_n\rightarrow f\left(\sigma^2\Sigma,\sigma^2\Sigma^{1/2}N_d\right)\quad\F_\infty\text{-mixing}
\end{displaymath}
and
\begin{displaymath}
D_n\rightarrow f\left(\sigma\Sigma^{1/2},\sigma^2\Sigma^{1/2}N_d\right)\quad\F_\infty\text{-mixing}
\end{displaymath}
as $n\to\infty$.
This concludes the proof.\hfill$\Box$
\medskip

\textit{Second proof of Theorem~\ref{Theorem:CLTsforAR}.} As explained above, the second proof of 
Theorem~\ref{Theorem:CLTsforAR} exploits stationarity and ergodicity of the autoregressive model. 
\smallskip

From now on we work with a bilateral i.i.d. sequence $\left(Z_n\right)_{n\in\ZZ}$
with $Z_1\in\LL^2\left(P\right)$, $E\left(Z_1\right)=0$ and $\sigma^2=\Var\left(Z_1\right)\in\left(0,\infty\right)$
(which represents a bilateral extension of the sequence of innovations $\left(Z_n\right)_{n\geq1}$). Let
$W_n=\left(Z_n,0,\ldots,0\right)^T,\,n\in\ZZ$.
\smallskip

Recall the dynamics
\begin{displaymath}
U_n=U_n\left(U_0\right)=B^nU_0+\sum_{j=1}^nB^{j-1}W_{n+1-j}=:B^nU_0+L_n\,,\quad n\geq0\,,
\end{displaymath}
established in \eqref{eq:Un}, where the summand $L_n$ is independent of the initial value $U_0$. The random variable 
$\U_0=\sum_{j=0}^\infty B^jW_{-j}$ is well-defined by the argument which proved \eqref{eq:sumfinite} and satisfies 
the integrability and independence 
assumptions on the initial values of the $U_n$-process. Write $\U_n=U_n\left(\U_0\right)$ for $n\geq0$ so that
$\left(\U_n\right)_{n\geq0}$ is the $U_n$-process for the initial value $\U_0$.
\smallskip

\textit{Step 1.}
The processes $\left(\U_n\right)_{n\geq0}$ and $\left(\U_{n-1}Z_n\right)_{n\geq1}$ are stationary and ergodic.
\smallskip

To see this we apply Lemma~\ref{App:preservation} with $\left(\X,\A\right)=\left(\RR,\B\left(\RR\right)\right)$,
$\left(\Y,\B\right)=\left(\RR^d,\B\left(\RR^d\right)\right)$ and $T_2=\NN_0$ and $\NN$, respectively. Note that for $n\geq0$
\begin{align*}
\U_n&=B^n\U_0+\sum_{j=1}^nB^{j-1}W_{n-\left(j-1\right)}=\sum_{j=0}^\infty B^{n+j}W_{-j}+\sum_{k=0}^{n-1}B^kW_{n-k}\\
&=\sum_{j=0}^\infty B^{n+j}W_{n-\left(n+j\right)}+\sum_{k=0}^{n-1}B^kW_{n-k}=
\sum_{k=n}^\infty B^kW_{n-k}+\sum_{k=0}^{n-1}B^kW_{n-k}\\
&=\sum_{k=0}^\infty B^kW_{n-k}
=\sum_{k=0}^\infty B^ke_1Z_{n-k}
\end{align*}
almost surely, where $e_1=\left(1,0,\ldots,0\right)^T\in\RR^d$ (moving average representation of $\U_n$). This suggests
the choice $\Lambda=\left\{x\in\RR^{\NN_0}:\sum_{j=0}^\infty\lVert B\rVert_B^j\left|x_j\right|<\infty\right\}$ and
$f:\RR^{\NN_0}\rightarrow\RR^d$,
\begin{displaymath}
f\left(x\right)=\left\{
\begin{array}{lll}
\displaystyle\sum_{j=0}^\infty B^je_1x_j&,&x\in\Lambda\\[20pt]
0                          &,&x\notin\Lambda\,.
\end{array}\right.
\end{displaymath}
Then clearly $\Lambda\in\B\left(\RR\right)^{\NN_0}$, $P\left(\left(Z_j\right)_{j\geq0}\in\Lambda\right)=
P\left(\left(Z_{n-j}\right)_{j\geq0}\in\Lambda\right)=1$ for every $n\in\ZZ$ and $f$ is measurable.
We obtain $\U_n=f\left(\left(Z_{n-j}\right)_{j\geq0}\right)$ almost surely for $n\geq0$. Since $\left(Z_n\right)_{n\in\ZZ}$
is stationary and ergodic (cf. \cite{Krengel}, Proposition~4.5 or \cite{Klenke}, Example~20.26), this
yields stationarity and 
ergodicity of $\left(\U_n\right)_{n\geq0}$ by Lemma~\ref{App:preservation}. Next choose $f_1:\RR^{\NN_0}\rightarrow\RR^d$,
$f_1\left(x\right)=x_0f\left(S^{\NN_0}\left(x\right)\right)$. Then
\begin{displaymath}
f_1\left(\left(Z_{n-j}\right)_{j\geq0}\right)=Z_nf\left(\left(Z_{n-1-j}\right)_{j\geq0}\right)=
Z_n\U_{n-1}
\end{displaymath}
almost surely for $n\geq1$ which again by Lemma~\ref{App:preservation} yields stationarity and ergodicity of
$\left(\U_{n-1}Z_n\right)_{n\geq1}$.
\smallskip

\textit{Step 2.} We also need the following asymptotic equivalences between the process $\U_n=U_n\left(\U_0\right),n\geq0$, and the process
$U_n=U_n\left(U_0\right),n\geq0$, with an arbitrary initial value $U_0$, namely
\begin{equation}
\label{eq:differenceUnUnbar}
\frac{1}{n}\sum_{j=1}^nU_{j-1}U_{j-1}^T-\frac{1}{n}\sum_{j=1}^n\U_{j-1}\U_{j-1}^T\rightarrow0\quad\text{almost surely}
\end{equation}
and
\begin{equation}
\label{eq:differenceUnZUnbarZ}
\frac{1}{\sqrt{n}}\sum_{j=1}^nU_{j-1}Z_j-\frac{1}{\sqrt{n}}\sum_{j=1}^n\U_{j-1}Z_j\rightarrow0\quad\text{almost surely}
\end{equation}
as $n\to\infty$.
\smallskip

As for \eqref{eq:differenceUnUnbar}, let $G_n=U_{n-1}U_{n-1}^T-\U_{n-1}\U_{n-1}^T$ for $n\geq1$. Then
\begin{align*}
G_n=&B^{n-1}\left(U_0U_0^T-\U_0\U_0^T\right)\left(B^{n-1}\right)^T+B^{n-1}\left(U_0-\U_0\right)L_{n-1}^T\\
&\qquad+L_{n-1}\left(U_0^T-\U_0^T\right)\left(B^{n-1}\right)^T
\end{align*}
so that
\begin{align*}
\lVert G_n\rVert_F\leq\lVert B^{n-1}\rVert_F^2\left(\lVert U_0\rVert^2+\lVert\U_0\rVert^2\right)
+2\lVert B^{n-1}\rVert_F\left(\lVert U_0\rVert+\lVert\U_0\rVert\right)\lVert L_{n-1}\rVert\,.
\end{align*}
Using
\begin{displaymath}
\lVert L_{n-1}\rVert\leq\sum_{j=1}^{n-1}\lVert B^{j-1}\rVert_F\left\lVert W_{n-j}\right\rVert=
\sum_{j=1}^{n-1}\lVert B^{j-1}\rVert_F\left\lvert Z_{n-j}\right\rvert\,,
\end{displaymath}
we find
\begin{align*}
E\left(\lVert L_{n-1}\rVert\right)\leq\sum_{j=1}^{n-1}\lVert B^{j-1}\rVert_FE\left(\left\lvert Z_{n-j}\right\rvert\right)
\leq\kappa_1E\left(\left\lvert Z_1\right\rvert\right)<\infty
\end{align*}
so that
\begin{displaymath}
\sum_{n=1}^\infty E\left(\lVert B^{n-1}\rVert_F\lVert L_{n-1}\rVert\right)=
\sum_{n=1}^\infty\lVert B^{n-1}\rVert_FE\left(\lVert L_{n-1}\rVert\right)\leq\kappa_1^2E\left(\left\lvert Z_1\right\rvert\right)<\infty\,.
\end{displaymath}
This yields $\lVert B^{n-1}\rVert_F\lVert L_{n-1}\rVert\rightarrow0$ almost surely. Combined with
$\left\lVert B^{n-1}\right\rVert_F^2\rightarrow0$ we obtain
$\lVert G_n\rVert_F\rightarrow0$ almost surely as $n\to\infty$. Consequently,
\begin{displaymath}
\left\lVert\frac{1}{n}\sum_{j=1}^nU_{j-1}U_{j-1}^T-\frac{1}{n}\sum_{j=1}^n\U_{j-1}\U_{j-1}^T\right\rVert_F
=\left\lVert\frac{1}{n}\sum_{j=1}^nG_j\right\rVert_F\leq\frac{1}{n}\sum_{j=1}^n\left\lVert G_j\right\rVert_F\rightarrow0
\end{displaymath}
almost surely as $n\to\infty$.
\smallskip

As for \eqref{eq:differenceUnZUnbarZ} we have for $n\geq1$
\begin{align*}
&\frac{1}{\sqrt{n}}\sum_{j=1}^nU_{j-1}Z_j-\frac{1}{\sqrt{n}}\sum_{j=1}^n\U_{j-1}Z_j=
\frac{1}{\sqrt{n}}\sum_{j=1}^n\left(U_{j-1}-\U_{j-1}\right)Z_j\\
&\qquad=\frac{1}{\sqrt{n}}\sum_{j=1}^nB^{j-1}\left(U_0-\U_0\right)Z_j
\end{align*}
so that
\begin{align*}
&\left\lVert\frac{1}{\sqrt{n}}\sum_{j=1}^nU_{j-1}Z_j-\frac{1}{\sqrt{n}}\sum_{j=1}^n\U_{j-1}Z_j\right\rVert\leq
\frac{1}{\sqrt{n}}\sum_{j=1}^n\lVert B^{j-1}\rVert_F\left(\lVert U_0\rVert+\lVert\U_0\rVert\right)\left\lvert Z_j\right\rvert\\
&\qquad\leq\frac{1}{\sqrt{n}}\left(\lVert U_0\rVert+\lVert\U_0\rVert\right)\sum_{j=1}^\infty
\lVert B^{j-1}\rVert_F\left\vert Z_j\right\rvert\rightarrow0
\end{align*}
almost surely as $n\to\infty$ in view of \eqref{eq:sumfinite}.
\smallskip

From Step~1 we know
that $\left(\U_{n-1}Z_n\right)_{n\geq1}$ is stationary and ergodic, and it is a martingale difference sequence w.r.t. the filtration
$\left(\G_n\right)_{n\geq0}$, where $\G_n=\sigma(U_0,Z_j,j\in\ZZ,j\leq n)$, with $E\left(\lVert\U_0Z_1\rVert^2\right)$ $=
E\left(\lVert\U_0\rVert^2\right)E\left(Z_1^2\right)<\infty$ 
by independence of $\U_0$ and $Z_1$ and $E\left(\lVert\U_0\rVert^2\right)<\infty$ and $E\left(Z_1^2\right)=\sigma^2<\infty$.
Therefore, by Corollary~\ref{Corollary:Stationarymartingaledifferences},
\begin{displaymath}
\frac{1}{\sqrt{n}}\sum_{k=1}^n\U_{k-1}Z_k\rightarrow E\left(\U_0Z_1\left(\U_0Z_1\right)^T\right)^{1/2}N_d\quad\G_\infty
\text{-mixing as }n\to\infty
\end{displaymath}
where $\G_\infty=\sigma\left(\bigcup_{n=0}^\infty\G_n\right)$,
with $E\left(\U_0Z_1\left(\U_0Z_1\right)^T\right)=E\left(Z_1^2\right)E\left(\U_0\U_0^T\right)=\sigma^4\Sigma$ so that
\begin{equation}
\label{eq:CLTinProof2}
\frac{1}{\sqrt{n}}\sum_{k=1}^n\U_{k-1}Z_k\rightarrow\sigma^2\Sigma^{1/2}N_d\quad\G_\infty
\text{-mixing as }n\to\infty\,.
\end{equation}
From \eqref{eq:CLTinProof2} and \eqref{eq:differenceUnZUnbarZ} we obtain for $U_n=U_n\left(U_0\right)$ with an
arbitrary initial value $U_0$ by an application of Theorem~3.18~(a) in \cite{HaeuslerLuschgy}
\begin{equation}
\label{eq:CLTforARwrongnorming2}
\frac{1}{\sqrt{n}}\sum_{k=1}^nU_{k-1}Z_k\rightarrow\sigma^2\Sigma^{1/2}N_d\quad\G_\infty
\text{-mixing and hence }\F_\infty\text{-mixing} 
\end{equation} 
as $n\to\infty$ since $\F_\infty\subset\G_\infty$.
Furthermore, Step~1 and the ergodic theorem (cf. Theorem~\ref{App:Birkhoff} with $X=\overline{U}$, $T=\NN_0$,
$f:\left(\RR^d\right)^{\NN_0}\rightarrow\RR^{d\times d}$, $f\left(x\right)=\pi_0\left(x\right)\pi_0\left(x\right)^T$,
where $\pi_0:\left(\RR^d\right)^{\NN_0}\rightarrow\RR^d$, $\pi_0\left(x\right)=x_0$ and
$\left\lVert f\left(\U\right)\right\rVert_F\leq\left\lVert\U_0\right\rVert^2\in\LL^1\left(P\right)$) imply
\begin{displaymath}
\frac{1}{n}\sum_{j=1}^n\U_{j-1}\U_{j-1}^T=\frac{1}{n}\sum_{j=0}^{n-1}\U_j\U_j^T
\rightarrow E\left(\U_0\U_0^T\right)=\sigma^2\Sigma\quad\text{almost surely as }
n\to\infty\,,
\end{displaymath}
so that by \eqref{eq:differenceUnUnbar}
\begin{equation}
\label{eq:UUTalmostsurely}
\frac{1}{n}\sum_{k=1}^nU_{k-1}U_{k-1}^T\rightarrow\sigma^2\Sigma\quad\text{almost surely as }n\to\infty\,.
\end{equation}
\smallskip

\textit{Step~3.} The remaining part of the proof is the same as Step~3 in the first proof of Theorem~\ref{Theorem:CLTsforAR},
because \eqref{eq:CLTforARwrongnorming2} is the same as \eqref{eq:CLTforARwrongnorming1} and \eqref{eq:UUTalmostsurely}
can be used instead of \eqref{eq:convergenceofUUT}.\hfill$\Box$
\medskip

From the equations 
\begin{displaymath}
\sum_{k=1}^nU_{k-1}U_{k-1}^T=\sum_{k=1}^mU_{k-1}U_{k-1}^T+\sum_{k=m+1}^nU_{k-1}U_{k-1}^T\quad\text{for all }n>m\geq1
\end{displaymath}
for symmetric and positive semidefinite matrices we see that $\sum_{k=1}^nU_{k-1}U_{k-1}^T$ is positive definite for all
$n\geq m$ whenever $\sum_{k=1}^mU_{k-1}U_{k-1}^T$ is positive definite. Therefore, the sequence $\left(\Omega_n\right)_{n\geq1}$
of events appearing in Step 3 the proof of Theorem~\ref{Theorem:CLTsforAR} is non-decreasing so that $P\left(\Omega_n\right)<1$ for
all $n\geq1$ or $P\left(\Omega_n\right)=1$ for all $n\geq n_0$ and some $n_0\geq1$. The following examples show that both
cases do indeed occur.
\medskip

\begin{Examples}
\label{Examples}
1. Let $P\left(U_0=0\right)>0$ and $P\left(Z_1=0\right)>0$. For all $n\geq1$
we have $\left\{Y_n=Y_{n-1}=\ldots=Y_1=Y_0=\ldots=Y_{-d+1}=0\right\}\subset\Omega_n^c$ so that
\begin{align*}
P\left(\Omega_n^c\right)&\geq P\left(Y_n=Y_{n-1}=\ldots=Y_1=Y_0=\ldots=Y_{-d+1}=0\right)\\
&\geq P\left(Z_n=\ldots=Z_1=0,U_0=0\right)=P\left(Z_1=0\right)^nP\left(U_0=0\right)>0\,.
\end{align*}
2. For $n\geq1$, let
\begin{displaymath}
\widetilde{Y}_n=
\begin{pmatrix}
Y_0     & Y_{-1}  & \ldots & Y_{-d+2}  & Y_{-d+1} \\
Y_1     & Y_0     & \ldots & Y_{-d+3}  & Y_{-d+2} \\    
Y_2     & Y_1     & \ldots & Y_{-d+4}  & Y_{-d+3} \\
\vdots  & \vdots  &        & \vdots    & \vdots   \\
Y_{d-1} & Y_{d-2} & \ldots & Y_1       & Y_0      \\
Y_d     & Y_{d-1} & \ldots & Y_2       & Y_1      \\
\vdots  & \vdots  &        & \vdots    & \vdots   \\
Y_{n-1} & Y_{n-2} & \ldots & Y_{n-d+1} & Y_{n-d}
 \end{pmatrix}
=\begin{pmatrix}
U_0^T  \\
U_1^T  \\
U_2^T  \\
\vdots \\
U_{d-1}^T \\
U_d^T \\
\vdots \\
U_{n-1}^T \\
\end{pmatrix}.
\end{displaymath}
Then $\widetilde{Y}_n^T\widetilde{Y}_n=\sum_{j=1}^nU_{j-1}U_{j-1}^T$ and 
$\text{rank}\left(\widetilde{Y}_n\right)=\text{rank}\left(\widetilde{Y}_n^T\widetilde{Y}_n\right)$.
We consider three cases. 
\smallskip

\textit{Case 1.} Assume $P\left(Y_0\neq0\right)=1$ and $P\left(Y_{-1}=\cdots=Y_{-d+1}=0\right)=1$.
Then with probability one $\widetilde{Y}_d$ is a lower triangular matrix with all diagonal elements different
from zero so that $\text{rank}\left(\widetilde{Y}_d\right)=d$ almost surely which implies $P\left(\Omega_d\right)=1$,
whence $P\left(\Omega_n\right)=1$ for all $n\geq d$.
\smallskip

\textit{Case 2.} Assume $P\left(U_0=0\right)=1$ and $P\left(Z_1\neq0\right)=1$. Then 
$Y_1=U_0^T\theta+Z_1=Z_1\neq0$ with probability one so that $\text{rank}\left(\widetilde{Y}_{d+1}\right)=d$ almost surely and
therefore $P\left(\Omega_{d+1}\right)=1$, whence $P\left(\Omega_n\right)=1$ for all $n\geq d+1$.\hfill$\Box$  

\textit{Case 3.} Assume that $P^{Z_1}$ is continuous. Then $\text{rank}\left(\widetilde{Y}_{2d}\right)=d$
almost surely. In fact, it follows that $P^{Y_n}$ is continuous for every $n\geq1$ and we have 
$d\geq\text{rank}\left(\widetilde{Y}_{2d}\right)\geq\text{rank}\left(D_d\right)$,
where $D_d$ denotes the $d\times d$-submatrix  of $\widetilde{Y}_{2d}$ consisting of the last $d$ rows of $\widetilde{Y}_{2d}$, i.e.
\begin{displaymath}
D_d=
\begin{pmatrix}
Y_d      & Y_{d-1}  & \ldots & Y_1     \\
\vdots   & \vdots   &        & \vdots  \\
Y_{2d-2} & Y_{2d-3} & \ldots & Y_{d-1} \\
Y_{2d-1} & Y_{2d-2} & \ldots & Y_d     \\
\end{pmatrix}.
\end{displaymath}
In order to show that $\text{rank}\left(D_d\right)=d$ almost surely consider for $1\leq j\leq d$ the
$j\times j$-submatrices $F_j$ of $D_d$ given by
\begin{displaymath}
F_j=
\begin{pmatrix}
Y_j      & Y_{j-1}  & \ldots & Y_1    \\
\vdots   & \vdots   &        & \vdots  \\
Y_{2j-2} & Y_{2j-3} & \ldots & Y_{j-1} \\
Y_{2j-1} & Y_{2j-2} & \ldots & Y_j     \\
\end{pmatrix},
\end{displaymath}
where $F_d=D_d$, and let $\Lambda_j=\left\{\det\left(F_j\right)\neq0\right\}$. Then by induction,
$P\left(\Lambda_j\right)=1$ for every $1\leq j\leq d$ since $\Lambda_1=\left\{Y_1\neq0\right\}$ and
$P\left(Y_1\neq0\right)=1$ and assuming $P\left(\Lambda_{j-1}\right)=1$ for some $2\leq j\leq d$, we obtain by
Laplace expansion of $\det\left(F_j\right)$  along the last row of $F_j$ that
\begin{align*}
\Lambda_j^c\cap\Lambda_{j-1}&=\left\{\det\left(F_j\right)=0\right\}\cap\Lambda_{j-1}\\
&=\left\{Y_{2j-1}\left(-1\right)^{j+1}\det\left(F_{j-1}\right)+f_1\left(Y_1,\ldots,Y_{2j-2}\right)=0\right\}
\cap\Lambda_{j-1}\\
&=\left\{Y_{2j-1}+f_2\left(Y_1,\ldots,Y_{2j-2}\right)=0\right\}\cap\Lambda_{j-1}\\
&=\left\{\vartheta_1Y_{2j-2}+\cdots+\vartheta_d Y_{2j-1-d}+Z_{2j-1}+f_2\left(Y_1,\ldots,Y_{2j-2}\right)=0\right\}\cap\Lambda_{j-1}\\
&=\left\{Z_{2j-1}=f_3\left(Y_m,\ldots,Y_{2j-2}\right)\right\}\cap\Lambda_{j-1}
\subset\left\{Z_{2j-1}=f_3\left(Y_m,\ldots,Y_{2j-2}\right)\right\}
\end{align*}
for $m=\left(2j-1-d\right)\wedge1$ and certain Borel measurable functions $f_1,f_2:\RR^{2j-2}\rightarrow\RR$ and
$f_3:\RR^s\rightarrow\RR$ with $s=2j-1-m=2j-2$ if $2j-2\geq d$ and $s=d$ if $2j-2 \leq d-1$ so that by independence of
$Z_{2j-1}$ and $\left(Y_m,\ldots,Y_{2j-2}\right)$,
\begin{align*}
P\left(\Lambda_j^c\right)&=P\left(\Lambda_j^c\cap\Lambda_{j-1}\right)\leq P\left(Z_{2j-1}=f_3\left(Y_m,\ldots,Y_{2j-2}\right)\right)\\
&=\int_{\RR^s}P\left(Z_{2j-1}=f_3\left(y_1,\ldots,y_s\right)\right)dP^{\left(Y_m,\ldots,Y_{2j-2}\right)}\left(y_1,\ldots,y_s\right)
=0\,.
\end{align*}
This implies $P\left(\Omega_{2d}^c\right)\leq P\left(\Lambda_d^c\right)=0$ and thus $P\left(\Omega_n\right)=1$ for all $n\geq2d$.
\hfill$\Box$
\end{Examples}
\bigskip

Finally, invertibility of the matrix $\Sigma$ is shown in the following lemma which already occurs in \cite{Anderson2}, Lemma~5.5.5.
\medskip

\begin{Lemma} 
\label{Lemma:Invertibility} 
The matrix $\sum_{j=0}^{d-1}B^j\widetilde{I}_d\left(B^j\right)^T$ is positive definite. Consequently, 
$\Sigma=\sum_{j=0}^{d-1}B^j\widetilde{I}_d\left(B^j\right)^T+\sum_{j=d}^\infty B^j\widetilde{I}_d\left(B^j\right)^T$
is positive definite as well and hence invertible.
\end{Lemma}
\medskip

\textit{Proof.} Since $\widetilde{I}_d$ is symmetric and idempotent,
$B^j\widetilde{I}_d\left(B^j\right)^T=B^j\widetilde{I}_d\left(B^j\widetilde{I}_d\right)^T$ so that 
$u^TB^j\widetilde{I}_d\left(B^j\right)^Tu=\left\lVert\left(B^j\widetilde{I}_d\right)^Tu\right\rVert^2$ for every $u\in\RR^d$.
Moreover, the first row of $\left(B^j\widetilde{I}_d\right)^T$ is given by the transposed first column of $B^j$ denoted by
$\left(B^j\right)_{{\scriptscriptstyle\bullet}1}$ while all other rows of $\left(B^j\widetilde{I}_d\right)^T$ are zero which implies
$\left\lVert\left(B^j\widetilde{I}_d\right)^Tu\right\rVert^2=\left\langle\left(B^j\right)_{{\scriptscriptstyle\bullet}1},u\right\rangle^2$. Consequently,
we obtain for every $u\in\RR^d$
\begin{displaymath}
u^T\left(\sum_{j=1}^{d-1}B^j\widetilde{I}_d\left(B^j\right)^T\right)u=
\sum_{j=1}^{d-1}u^TB^j\widetilde{I}_d\left(B^j\right)^Tu=
\sum_{j=1}^{d-1}\left\langle\left(B^j\right)_{{\scriptscriptstyle\bullet}1},u\right\rangle^2\,.
\end{displaymath}
For every $j=0,\cdots,d-1$ the first column of the $d\times d$-matrix $B^j$ satisfies $\left(B^j\right)_{j+1,1}=1$ and
$\left(B^j\right)_{j+k,1}=0$ for $k=2,\cdots,d-j$. This follows by induction from 
$\left(B^0\right)_{{\scriptscriptstyle\bullet}1}=\left(I_d\right)_{{\scriptscriptstyle\bullet}1}=\left(1,0,\ldots,0\right)^T$ and
\begin{displaymath}
\left(B^{j+1}\right)_{{\scriptscriptstyle\bullet}1}=\left(BB^j\right)_{{\scriptscriptstyle\bullet}1}=
\left(\left\langle\theta,\left(B^j\right)_{{\scriptscriptstyle\bullet}1}\right\rangle,\left(B^j\right)_{11},\ldots,\left(B^j\right)_{d-1,1}\right)^T\,.
\end{displaymath}
Now assume $\sum_{j=0}^{d-1}\left\langle\left(B^j\right)_{{\scriptscriptstyle\bullet}1},u\right\rangle^2=0$ for some $u\in\RR^d$. Hence 
$\left\langle\left(B^j\right)_{{\scriptscriptstyle\bullet}1},u\right\rangle=0$ for every $j=0,\ldots,d-1$ or equivalently $Du=0$, where the
$j$-th row of the $d\times d$-matrix $D$ denoted by $D_{j{\scriptscriptstyle\bullet}}$ is given by 
$D_{j{\scriptscriptstyle\bullet}}=\left(\left(B^{j-1}\right)_{{\scriptscriptstyle\bullet}1}\right)^T$ for $j=1,\ldots,d$. 
Then $D$ is a lower triangular matrix 
with all diagonal elements equal to 1 so that $\det\left(D\right)=1\neq0$. We conclude that $u=0$ is the only
solution of the equation $Du=0$.\hfill$\Box$
\bigskip\smallskip

{\rm\bf\large Appendix}
\bigskip

\setcounter{Theorem}{0}

\setcounter{section}{1}

\renewcommand{\thesection}{\Alph{section}}

We recall some basic facts about stationarity and ergodicity of random processes. Let $\left(\X,\A\right)$
denote a measurable space and let $T\in\left\{\NN,\NN_0,\ZZ\right\}$. The shift $S=S^T$ on $\X^T$ is defined by
$S\left(x\right)=S\left(\left(x_n\right)_{n\in T}\right)=\left(x_{n+1}\right)_{n\in T}$. It is
$\left(\A^T,\A^T\right)$-measurable. Let $\A^T\left(S\right)=\left\{A\in\A^T:S^{-1}\left(A\right)=A\right\}$ denote the
$\sigma$-field of invariant (more precisely, $S$-invariant) measurable subsets of $\X^T$. Let $X=\left(X_n\right)_{n\in T}$
be an $\X$-valued process defined on some probability space $\left(\Omega,\F,P\right)$ 
and note that $X$ can be seen as an $\X^T$-valued random vector. The process $X$ is called
(\textit{strictly}) \textit{stationary} if $P^{S\left(X\right)}=P^X$ and \textit{ergodic} if
$P^X\left(\A^T\left(S\right)\right)=\left\{0,1\right\}$.  
\smallskip 

In the sequel let $\left(\Y,\B\right)$ be a further measurable space and let $T_1,T_2\in\left\{\NN,\NN_0,\ZZ\right\}$.

\begin{Lemma}
\label{App:equivariance}
Let $g:\X^{T_1}\rightarrow\Y^{T_2}$ be an $\left(\A^{T_1},\B^{T_2}\right)$-measurable map which is equivariant in the sense of
$g\circ S^{T_1}=S^{T_2}\circ g$. Then, if $X=\left(X_n\right)_{n\in T_1}$ is an $\X$-valued stationary process, the
$\Y$-valued process $Y:=g\left(X\right)$ is stationary. If $X$ is ergodic, then $Y$ is ergodic. 
\end{Lemma}
\medskip

Here the equation
$Y=g\left(X\right)$ means
\begin{displaymath}
Y_n=\pi_n^{T_2}\left(Y\right)=\pi_n^{T_2}\left(g\left(X\right)\right)=\left(\pi_n^{T_2}\circ g\right)\left(X\right)
=:g_n\left(X\right)=g_n\left(\left(X_n\right)_{n\in T_1}\right)
\end{displaymath}
for all $n\in T_2$, where $\pi_n^{T_2}:\Y^{T_2}\rightarrow\Y$ denote the projections on $\Y^{T_2}$.
\medskip

\textit{Proof.} (see \cite{Kallenberg}, Lemma A1.1) Assume that $X$ is stationary. Then
\begin{displaymath}
P^{S^{T_2}\left(Y\right)}=P^{S^{T_2}\circ g\left(X\right)}=P^{g\circ S^{T_1}\left(X\right)}
=\left(P^{S^{T_1}\left(X\right)}\right)^g=\left(P^X\right)^g=P^{g\left(X\right)}=P^Y
\end{displaymath}
so that $Y$ is stationary. Next we note that 
$g^{-1}\left(\B^{T_2}\left(S^{T_2}\right)\right)\subset\A^{T_1}\left(S^{T_1}\right)$. In fact, if
$H\in\B^{T_2}\left(S^{T_2}\right)$, then
\begin{align*}
\left(S^{T_1}\right)^{-1}\left(g^{-1}\left(H\right)\right)&=\left(g\circ S^{T_1}\right)^{-1}\left(H\right)=
\left(S^{T_2}\circ g\right)^{-1}\left(H\right)=g^{-1}\left(\left(S^{T_2}\right)^{-1}\left(H\right)\right)\\
&=g^{-1}\left(H\right)\,.
\end{align*}
Hence, if $X$ is ergodic, then
\begin{displaymath}
P^Y\left(\B^{T_2}\left(S^{T_2}\right)\right)=P^X\left(g^{-1}\left(\B^{T_2}\left(S^{T_2}\right)\right)\right)\subset
P^X\left(\A^{T_1}\left(S^{T_1}\right)\right)=\left\{0,1\right\}
\end{displaymath}
which shows that $Y$ is ergodic.\hfill$\Box$ 
\medskip

Equivariance of $g=\left(g_n\right)_{n\in T_2}$ may be read as follows.
\medskip

\begin{Lemma}
\label{App:equivalence}
Assume $T_1=\ZZ$ if $T_2=\ZZ$. Let $g:\X^{T_1}\rightarrow\Y^{T_2}$ be measurable and let $S=S^{T_1}$.

{\rm (a)} Let $T_2\in\left\{\NN_0,\ZZ\right\}$. Then $g$ is equivariant if and only if $g_n=h\circ S^n$ for every
$n\in T_2$ and some measurable map $h:\X^{T_1}\rightarrow\Y$. The same holds true in case $T_2=\NN$ and
$T_1=\ZZ$.

{\rm (b)} Let $T_2=\NN$ and $T_1\in\left\{\NN,\NN_0\right\}$. Then g is equivariant if and only if 
$g_n=h\circ S^{n-1}$ for every $n\in\NN$ and some measurable map $h:\X^{T_1}\rightarrow\Y$. A
sufficient condition for equivariance of $g$ is $g_n=h\circ S^n$ for every $n\in\NN$ and some
measurable map $h$.
\end{Lemma}
\medskip

\textit{Proof.} Note that equivariance of $g$ means that $g_n\circ S=g_{n+1}$ for every $n\in T_2$. In any case,
the condition $g_n=h\circ S^n$ for every $n\in T_2$ implies equivariance of $g$. In fact, we obtain
\begin{displaymath}
g_{n+1}=h\circ S^{n+1}=h\circ S^n\circ S=g_n\circ S
\end{displaymath}
for every $n\in T_2$ so that $g$ is equivariant. Also the condition $g_n=h\circ S^{n-1}$ for every $n\in T_2=\NN$ 
implies equivariance of $g$ since 
\begin{displaymath}
g_{n+1}=h\circ S^n=h\circ S^{n-1}\circ S=g_n\circ S
\end{displaymath}
for every $n\in T_2$.
\smallskip

(a) Assume that $g$ is equivariant. Let $T_2\in\left\{\NN_0,\ZZ\right\}$. Then we obtain from
$g_{n+1}=g_n\circ S,\,n\in T_2$ by induction $g_n=g_0\circ S^n,\,n\in\NN_0$. If $T_2=\ZZ$, then $T_1=\ZZ$ and
by induction $g_{-n}\circ S^n=g_0,\,n\in\NN_0$ so that $g_{-n}=g_0\circ S^{-n}$ for every $n\in\NN_0$ since
$S=S^\ZZ$ is bijective. If $T_2=\NN$ and $T_1=\ZZ$, then by induction $g_n=g_1\circ S^{n-1},\,n\in T_2=\NN$.
Choosing $h=g_1\circ S^{-1}$ yields $g_n=h\circ S^n$ for every $n\in T_2=\NN$.
\smallskip

(b) If $g$ is equivariant, then by induction $g_n=g_1\circ S^{n-1},\,n\in\NN$.\hfill$\Box$
\medskip

As a consequence we obtain the following preservation of stationarity and ergodicity.
\medskip

\begin{Lemma}
\label{App:preservation}
Let $T_1=\ZZ$ and let $X=\left(X_n\right)_{n\in\ZZ}$ be an $\X$-valued stationary and ergodic process. Let
$f:\X^{\NN_0}\rightarrow\Y$ be an $\left(\A^{\NN_0},\B\right)$-measurable map and let $Y_n=f\left(\left(X_{n-j}\right)_{j\geq0}\right)$
for every $n\in T_2$. Then $Y=\left(Y_n\right)_{n\in T_2}$ is stationary and ergodic.
\end{Lemma}
\medskip

\textit{Proof.} The map $\varphi:\X^\ZZ\rightarrow\X^{\NN_0},\,\varphi\left(x\right)=\left(x_{-j}\right)_{j\geq0}$ is
$\left(\A^\ZZ,\A^{\NN_0}\right)$-measurable since $\pi_j^{\NN_0}\left(\varphi\left(x\right)\right)=x_{-j}=\pi_{-j}^\ZZ\left(x\right)$
for $j\geq0$, where $\pi_j^{\NN_0}$ and $\pi_{-j}^\ZZ$ are projections on $\X^{\NN_0}$ and $\X^\ZZ$, respectively, and
$\pi_{-j}^\ZZ$ is $\left(\A^\ZZ,\A\right)$-measurable. Define $g=\left(g_n\right)_{n\in T_2}:\X^\ZZ\rightarrow\Y^{T_2}$ by 
$g_n=f\circ\varphi\circ\left(S^\ZZ\right)^n$. Then $g$ is equivariant by Lemma~\ref{App:equivalence} and we have
$g_n\left(x\right)=f\left(\left(x_{n-j}\right)_{j\geq 0}\right)$ for every $x\in\X^\ZZ,\,n\in T_2$ so that $Y=g\left(X\right)$.
The assertion follows from Lemma~\ref{App:equivariance}.\hfill$\Box$
\medskip

The pointwise ergodic theorem of Birkhoff in this setting reads as follows. Let $T\in\left\{\NN,\NN_0,\ZZ\right\}$
and let $\I:=\A^T\left(S\right)$, where $S=S^T$.
\medskip

\begin{Theorem}
\label{App:Birkhoff}
{\rm(}Ergodic theorem, Birkhoff\;{\rm)} Let $X=\left(X_n\right)_{n\in T}$ be an $\X$-valued stationary process. If
$f\in\LL^1\left(P^X\right)$, then
\begin{displaymath}
\frac{1}{n}\sum_{j=0}^{n-1}f\left(\left(X_{k+j}\right)_{k\in T}\right)=
\frac{1}{n}\sum_{j=0}^{n-1}f\circ S^j\left(X\right)\rightarrow
E_{P^X}\left(f|\I\right)\circ X=E\left(f\left(X\right)|\I_X\right)
\end{displaymath}
almost surely as $n\to\infty$ with $\I_X=X^{-1}\left(\I\right)$. The same holds true for measurable maps $f:\X^T\rightarrow\RR^{k\times m}$
with $\left\lVert f\left(X\right)\right\rVert_F\in\LL^1\left(P\right)$ by componentwise application of the result for
real-valued functions $f$.
\end{Theorem}
\medskip

\textit{Proof.} See \cite{Krengel}, Theorem 2.3, or \cite{Klenke}, Theorem~20.14.\hfill$\Box$
\medskip

Note that by Lemmas~\ref{App:equivariance} and \ref{App:equivalence}, the processes $Y_j=f\circ S^j\left(X\right),\,j\geq0$ and
$Y_j=f\circ S^{j-1}\left(X\right),\,j\geq1$ are stationary.
\bigskip

\renewcommand\refname

\end{document}